\documentclass[aop,MSNbibl,seceqn,dvips]{arximspdf}
\usepackage{accents}

\doi{10.1214/13-AOP848} %kopijuoti is PTS
\volume{42}
\issue{5}
\pubyear{2014}
\firstpage{2161}
\lastpage{2196}

\makeatletter

\newcommand{\nlimsup}{\mathop{\operatorname{\overline{lim}}}}

\newtheorem{theorem}{Theorem}[section]
\newtheorem{lemma}[theorem]{Lemma}

\newtheorem{corollary}[theorem]{Corollary}

\newproclaim{assumption}{Assumption}[section]
\newproclaim{definition}{Definition}[section]
\newproclaim{remark}{Remark}[section]

\newproclaim{example}{Example}[section]

\newcommand{\tr}{\operatorname{tr}}

\newcommand\V{\check{\hspace*{-4pt}\phantom{x}}}

\newcommand\bbeta{\bolds{\beta}}

\newcommand\bR{\mathbb{R}}

\newcommand\cF{\mathcal{F}}
\newcommand\cO{\mathcal{O}}
\newcommand\cP{\mathcal{P}}
\newcommand\cSk{\mathcal{S}k}

\newcommand\frA{\mathfrak{A}}
\newcommand\frB{\mathfrak{B}}

\newcommand\Bsigma{\bolds{\sigma}}
\newcommand\Ba{\mathbf{a}}
\newcommand\Bb{\mathbf{b}}
\newcommand\Bc{\mathbf{c}}
\newcommand\Bf{\mathbf{f}}
\newcommand\BL{\mathbf{L}}

\newcommand\dist{\operatorname{dist}}

\newcommand\infsup{\mathop{\operatorname{inf\,sup}}}
\newcommand\supinf{\mathop{\operatorname{sup\,inf}}}
\newcommand{\eqref}[1]{(\ref{#1})}
\makeatother

\begin{document}
\begin{frontmatter}

\title{On regularity properties and approximations of~value functions
for stochastic differential~games~in~domains}
\runtitle{Smoothness of value functions}

\begin{aug}
\author{\fnms{N.~V.} \snm{Krylov}\corref{}\thanksref{t1}\ead[label=e1]{krylov@math.umn.edu}}
\thankstext{t1}{Supported in part by NSF Grant DMS-1160569.}
\runauthor{N.~V. Krylov}
\affiliation{University of Minnesota}
\address{Department of Mathematics\\
University of Minnesota\\
127 Vincent Hall\\
Minneapolis, Minnesota 55455\\
USA\\
\printead{e1}}
\end{aug}

% HISTORY:
\received{\smonth{8} \syear{2012}}
\revised{\smonth{3} \syear{2013}}

% ABSTRACT
%
\begin{abstract}
We prove that for any
constant $K\geq1$,
the value functions
for time homogeneous
stochastic differential games in the whole space
can be approximated
up to a constant over $K$ by value functions
whose second-order derivatives are \emph{bounded} by a
constant times $K$.

On the way of proving this result
we prove that the value functions for
stochastic differential games in domains and in the whole
space admit estimates
of their Lipschitz constants in a variety of
settings.
\end{abstract}

% KEYWORDS
% Pirmas kwd is didziosios raides
%

\end{frontmatter}

%s1 #&#
%s1 ###
\section{Introduction}\label{sec1}

In this paper we prove that for any
constant $K\geq1$,
the value functions
for time homogeneous
stochastic differential games in the whole space
can be approximated
up to a constant over $K$ by the value functions
whose second-order derivatives are \emph{bounded} by a
constant times $K$ (see Theorem~\ref{theorem 6.3.1} and Remark~\ref{remark 6.3.1}).
To prove Theorem~\ref{theorem 6.3.1}
we needed a few auxiliary facts organized
in~\cite{Kr_1} and \cite{Kr_2}, so that the goal to prove
this theorem
was the major driving force of the series of three articles.
Along the way some fruitful ideas were developed,
leading,
in particular, first to understanding from probabilistic point of
view and then to proving in purely
PDE terms the fact that one can find in
$C^{1+\alpha}$ viscosity solutions
of the uniformly nondegenerate Isaacs parabolic equations
with coefficients measurable in time and
VMO in $x$; see \cite{Kr_3}.
It would be extremely interesting to find a proof
of this fact based on the theory of viscosity solutions
in the situation of discontinuous coefficients, although
in the case of continuous ones such a proof
was given by \'Swi{\c e}ch \cite{Sw97}.

In terms of the corresponding Isaacs
equations the approximation
in Theorem~\ref{theorem 6.3.1}
is done in such a way that the equations are modified
only for large values of the derivatives of the value
functions. Such approximation of stochastic games
can be useful while evaluating the value functions
numerically because one can expect that approximations
might be more accurate if the approximating function
is more regular.

Two main tools are used. One is the stochastic dynamic
principle with randomized stopping times, and another
is based on estimates of the Lipschitz constants
of the value functions.

The dynamic programming principle
we use is proved in \cite{Kr_2} and originated in the work
by Fleming and Souganidis \cite{FS89}; see also
Kovats \cite{Ko09} and \'Swi{\c e}ch \cite{Sw96}.

Here we concentrate on proving the Lipschitz
continuity of the value functions for
time homogeneous stochastic differential games
in domains and in the whole space and on proving
the above mentioned approximation result, which
is a particular case of a conjecture from \cite{Kr12.2}.

There is an enormous literature treating smoothness
properties for \emph{controlled} diffusion processes or, from
analytical point of view, for fully nonlinear equations under
convexity assumptions. We are going to focus only on
stochastic differential games for which there is not much
known concerning the regularity of the value function
in more or less general case.

Ishii and Lions in \cite{IL90} prove the
Lipschitz continuity for viscosity solutions
of fully nonlinear uniformly nondegenerate equations.
Earlier
Trudinger in \cite{Tr89} proved that
the first derivatives are, actually, H\"older
continuous. The same result under somewhat more
restrictive assumptions can be found in the book
\cite{CC95} by Caffarelli and Cabr\'e.
Further results on Lipschitz continuity, still
for uniformly nondegenerate case,
are contained and referred to in
\'{S}wi\c{e}ch \cite{Sw97}, Vitolo \cite{Vi10} and Krylov~\cite{Kr_3}.

We deal with global and local estimates only
for the Isaacs equations
in contrast with the more general
equations in the above mentioned references, which reduce
to the Isaacs equations only if the equation
is determined by the so-called
boundedly inhomogeneous functions. Our methods are
also different from the methods of the above cited
articles where the authors rely on the theory of viscosity
solutions. Our solutions are given as value
functions of stochastic differential games,
and we use probabilistic methods, with the main tool
being based on different probabilistic representations
for the value functions at different points.
This is very close to using the so-called
quasiderivatives of solutions of stochastic
equations in the theory of controlled diffusion
processes, which can be traced down starting from \cite{Kr89}.
We could also use quasiderivatives in this article,
but this would require more work, and what we are
actually using can be called the method of
quasidifferences. In the author's opinion the methods
of this article can be also applied to proving
interior first derivatives estimates for
degenerate equations similar to those
in \cite{Zh12}
when the boundary data are only
Lipschitz continuous and processes are not
uniformly nondegenerate.
Just in case, observe that there are no global
gradient estimates even for the equation $\Delta u=0$
in a ball if the boundary data are only Lipschitz continuous.

Even though our stochastic differential games
are assumed to be uniformly nondegenerate, one
of our main results, Theorem~\ref{theorem 5.22.1}, is about
estimates of the Lipschitz constant \emph{independent}
of the constant of nondegeneracy.
The author is not aware
of any analytical proof of it.
The only results similar to the one mentioned above that
the author is aware of are contained in Barles \cite{Ba91}.
We discuss them in detail in Remark~\ref{remark 5.31.4}.

We also prove two estimates which do depend on the constant of
nondegeneracy:
one is global, Theorem~\ref{theorem 5.19.1}, and
another is local, Theorem~\ref{theorem 6.4.1}.
These results are much weaker than the ones in
\cite{Tr89}. The emphasis here is to show that
probabilistic methods can use nondegeneracy
in an efficient way. Of course, Theorem~\ref{theorem 5.22.1}
contains Theorem~\ref{theorem 5.19.1}, the proof of the latter
is given just because it is short, instructive and requires less
machinery.

The main results of the paper are stated in Section~\ref{section 2.26.3}. Section~\ref{section 5.29.6}
contains their discussion continued in Section~\ref{section 3.8.1} where we describe some ideas
behind our arguments.
In Section~\ref{section 5.20.1}
we show that the value function admits many
representations. In Section~\ref{section 5.20.2} we prove auxiliary results
aimed at estimating the difference of value
function at close points when different
probabilistic representations are taken
for those points. The result of
Section~\ref{section 5.20.1}, in a very rough
form, is used in Section~\ref{section 5.29.1} to prove
Theorem~\ref{theorem 5.19.1}. In
Section~\ref{section 6.4.1} we prove Theorem~\ref{theorem 6.4.1} about interior estimates.
A very short Section~\ref{section 6.4.2} contains the
proof of Theorem~\ref{theorem 5.22.1} about estimates
independent of the constant of nondegeneracy.
It is short because the main ideas have already been
given in Section~\ref{section 5.20.1}.
In the final, and again short, Section~\ref{section 6.8.1}
we prove Theorem~\ref{theorem 6.3.1}.

The author is very grateful to the referees for their comments
which helped improve the presentation of this paper.

%s2 #&#
%s2 ###
\section{Main results}
\label{section 2.26.3}

Let $\bR^{d}=\{x=(x_{1},\ldots,x_{d})\}$
be a $d$-dimensional Euclidean space, and let $d_{1}\geq d$
be an integer. Denote by $\cO$ the set of
$d_{1}\times d_{1}$ orthogonal
matrices, fix an integer $k\geq1$ and
assume that we are given separable metric spaces
$A$ and $B$ and let,
for each $\alpha\in A$, $\beta\in B$ and
$p\in\bR^{k}$, the following
functions on $\bR^{k}\times\bR^{d}$ be given:
\begin{longlist}[(iii)]
\item[(i)] $d\times d_{1}$
matrix-valued $ \sigma^{\alpha\beta}(p,x) =
(\sigma^{\alpha\beta}_{ij}(p,x))$;

\item[(ii)] $\cO$-valued
function $P^{\alpha\beta}(x,y)$, $\bR^{k}$-valued
function $p^{\alpha\beta}(x,y)$ and real-valued
function $r^{\alpha\beta}(x,y)$;

\item[(iii)]
$\bR^{d}$-valued $ b^{\alpha\beta}(p,x) =
( b^{\alpha\beta}_{i }(p,x))$;

\item[(iv)]
real-valued functions
$ c^{\alpha\beta}(p,x)\geq0$,
$ f^{\alpha\beta}(p,x)$ and
$g(x)$.
\end{longlist}

Define
\[
a^{\alpha\beta}(p,x):=(1/2) \sigma^{\alpha\beta}(p,x) \bigl(
\sigma^{\alpha\beta}(p,x) \bigr)^{*}.
\]
Also set
\[
(\sigma,a,b,c,f)^{\alpha\beta}(x) =(\sigma,a,b,c,f)^{\alpha\beta}(0,x),
\]
and note that for our first main result, Theorem~\ref{theorem 5.19.1},
only these values of
$\sigma,a,b,c,f$ are relevant,
and the parameters $r,p,P$ are not present. These
parameters
are important in Theorem~\ref{theorem 5.22.1}.
The role of these parameters
is discussed in Remark~\ref{remark 5.29.2}
and Example~\ref{example 3.9.1}
concerning $P$, in
Remarks \ref{remark 5.21.2}, \ref{remark 6.1.1}
and Example~\ref{example 3.9.2} concerning $r$
and in
Remark~\ref{remark 6.1.1} concerning $p$.

Fix some constants $K_{0},K_{1}\in[0,\infty)$,
and $\delta_{0} \in(0,1]$.

%
%as2.1 #&#
\begin{assumption}
\label{assumption 5.19.1}
(i) The functions
$(\sigma,a,b,c,f)^{\alpha\beta}(p,x)$
and $p^{\alpha\beta}(x,y)$ are continuous with respect to
$\beta\in B$ for each $(\alpha,p,x,y)$ and continuous with respect
to $\alpha\in A$ uniformly with respect to $\beta\in B$
for each $(p,x,y)$. Furthermore, they are Borel measurable functions
of $(p,x,y)$ for each
$(\alpha,\beta)$ and they are bounded by~$K_{0}$.\vspace*{-6pt}

\begin{longlist}[(iii)]
\item[(ii)] The functions $r^{\alpha\beta}(x,y)$ and
$P^{\alpha\beta}(x,y)$
are bounded by constant $K_{0}$, they are Borel measurable
with respect to all variables, and along with $p^{\alpha\beta}(x,y)$
they are
Lipschitz
continuous with respect to $x$ with Lipschitz constant $K_{1}$, and
\[
r^{\alpha\beta}(x,x)\equiv1,\qquad  p^{\alpha\beta}(x,x)\equiv0,\qquad  P^{\alpha\beta}(x,x)
\equiv I,
\]
where $I$ is the $d_{1}\times d_{1}$-identity matrix.
The function $p^{\alpha\beta}(x,y)$ is uniformly continuous
with respect to $y$ uniformly with respect to $(\alpha,\beta,x)$.

\item[(iii)] The functions $\sigma^{\alpha\beta}(p,x)$,
$b^{\alpha\beta}(p,x)$, $ c^{\alpha\beta}(p,x)$ and
$ f^{\alpha\beta}(p,x)$ are Lipschitz
continuous with respect to $(p,x)$ with Lipschitz
constant $K_{1}$. We have $\|g\|_{C^{2}(\bR^{d})}\leq K_{1}$.

\item[(iv)] For
any
$\alpha\in A$, $\beta\in B$, $x,\lambda\in\bR^{d}$ and
$p\in\bR^{k}$, we have
\[
a^{\alpha\beta}_{ij}(p,x)\lambda_{i}
\lambda_{j}\geq\delta_{0}|\lambda|^{2}.
\]
\end{longlist}
\end{assumption}
The reader understands, of course, that the summation
convention is adopted throughout the article.

Let $(\Omega,\cF,P)$ be a complete probability space,
let $\{\cF_{t},t\geq0\}$ be an increasing filtration
of $\sigma$-fields $\cF_{t}\subset\cF$ such that
each $\cF_{t}$ is complete with respect to $\cF,P$ and let
$w_{t},t\geq0$, be a standard $d_{1}$-dimensional Wiener process
given on $\Omega$ such that $w_{t}$ is a Wiener process
relative to the filtration $\{\cF_{t},t\geq0\}$.

The set of progressively measurable $A$-valued
processes $\alpha_{t}=\alpha_{t}(\omega)$ is denoted by $\frA$.
Similarly we define $\frB$
as the set of $B$-valued progressively measurable functions.
By $ \mathbb{B}$ we denote
the set of $\frB$-valued functions
$ \bbeta(\alpha_{\cdot})$ on $\frA$
such that, for any $T\in(0,\infty)$ and any $\alpha^{1}_{\cdot},
\alpha^{2}_{\cdot}\in\frA$ satisfying
%
%
%e2.1 #&#
%e2.1 ###
\begin{equation}
\label{4.5.4} P \bigl( \alpha^{1}_{t}=
\alpha^{2}_{t} \mbox{ for almost all } t\leq T \bigr)=1,
\end{equation}
we have
\[
P \bigl( \bbeta_{t} \bigl(\alpha^{1}_{\cdot}
\bigr)= \bbeta_{t} \bigl(\alpha^{2}_{\cdot} \bigr)
\mbox{ for almost all } t\leq T \bigr)=1.
\]

For $\alpha_{\cdot}\in\frA$, $\beta_{\cdot}\in\frB$ and
$x\in\bR^{d}$ introduce
$x^{\alpha_{\cdot}\beta_{\cdot} x}_{t}$
as a unique solution of the It\^o equation
%
%
%e2.2 #&#
%e2.2 ###
\begin{equation}
\label{5.11.1} x_{t}=x+\int_{0}^{t}
\sigma^{\alpha_{s}
\beta_{s} }( x_{s}) \,dw_{s} +\int
_{0}^{t}b^{\alpha_{s}
\beta_{s} }( x_{s}) \,ds,
\end{equation}
and denote
\[
\phi^{\alpha_{\cdot}\beta_{\cdot} x}_{t} =\int_{0}^{t}c^{\alpha_{s}
\beta_{s} }
\bigl( x^{\alpha_{\cdot}
\beta_{\cdot} x}_{s} \bigr) \,ds.
\]
Next, fix a domain $D\subset\bR^{d}$,
define $\tau^{\alpha_{\cdot}\beta_{\cdot} x}$ as the first exit
time of $x^{\alpha_{\cdot}
\beta_{\cdot} x}_{t}$ from $D$
($\tau^{\alpha_{\cdot}\beta_{\cdot} x}=\infty$
if $D=\bR^{d}$) and introduce
%
%
%e2.3 #&#
%e2.3 ###
\begin{equation}
\label{2.12.2} v(x)=\infsup_{\bbeta\in\mathbb{B}\,  \alpha_{\cdot}\in\frA} E_{x}^{\alpha_{\cdot}\bbeta(\alpha_{\cdot})}
\biggl[\int_{0}^{\tau} f( x_{t})e^{-\phi_{t}}
\,dt+g(x_{\tau})e^{-\phi_{\tau}} \biggr],
\end{equation}
where the indices $\alpha_{\cdot}$, $\bbeta$, and $x$
at the expectation sign are written to mean that
they should be placed inside the expectation sign
wherever and as appropriate, that is,
\begin{eqnarray*}
&&E_{x}^{\alpha_{\cdot}\beta_{\cdot}} \biggl[\int_{0}^{\tau}
f( x_{t})e^{-\phi_{t}} \,dt+g(x_{\tau})e^{-\phi_{\tau}}
\biggr]
\\
&&\qquad:= E \biggl[ g \bigl(x^{\alpha_{\cdot}\beta_{\cdot} x} _{\tau^{\alpha_{\cdot}\beta_{\cdot} x}} \bigr) e^{-\phi^{\alpha_{\cdot}\beta_{\cdot} x}
_{\tau^{\alpha_{\cdot}\beta_{\cdot} x}}} +
\int_{0}^{\tau^{\alpha_{\cdot}\beta_{\cdot} x}} f^{\alpha_{t}\beta_{t} } \bigl(
x^{\alpha_{\cdot}\beta_{\cdot} x}_{t} \bigr) e^{-\phi^{\alpha_{\cdot}\beta_{\cdot} x}_{t}} \,dt \biggr].
\end{eqnarray*}
Observe that $v(x)=g(x)$ in $\bR^{d}\setminus D$.
Next, introduce
\[
L^{\alpha\beta}u(p,x)=a^{\alpha\beta}_{ij}(p,x) D_{ij}u(x)
+b^{\alpha\beta}_{i}(p,x)D_{i}u(x)-c^{\alpha\beta}(p,x)
u(x),
\]
where $D_{i}=\partial/(\partial x_{i})$, $D_{ij}=D_{i}D_{j}$
and note for orientation that
$v$ is a viscosity solution
of the corresponding Isaacs equation
\[
\supinf_{\alpha\in A\,  \beta\in B} \bigl[ L^{\alpha\beta}u(0,x)+f^{\alpha\beta}(x)
\bigr]=0, \qquad x\in D.
\]
This fact which will not play any role here
is proved
in \cite{Ko09} for bounded domains.

Our first main result is the following.

%
%th2.1 #&#
\begin{theorem}
\label{theorem 5.19.1}
Under the above assumptions also suppose that
either $D$ is bounded and
satisfies the uniform exterior ball condition, or
$D=\bR^{d}$ and there is a constant $\delta_{1}>0$ such that
$c^{\alpha\beta}(x)\geq\delta_{1}$.

Then
$v$ is Lipschitz
continuous in $\bR^{d}$ with Lipschitz constant
depending only on~$D$, $K_{0},K_{1}$, $\delta_{0}$
and $\delta_{1}$.
\end{theorem}

The above setting and notation follow \cite{Kr_2}
and, as there, we convince ourselves that the definition of
$v$ makes sense, and $v$ is bounded.

Here is a result about interior smoothness of $v$.
%
%
%th2.2 #&#
\begin{theorem}
\label{theorem 6.4.1}
Let $D$ be bounded and in Assumption~\ref{assumption 5.19.1}\textup{(iii)}
replace the requirement $\|g\|_{C^{2}(\bR^{d})}\leq K_{1}$
with the requirement that $g$ is continuous. Then
$v$ is Lipschitz continuous
on any compact set $\Gamma\subset D$.
\end{theorem}

As we have pointed out in the \hyperref[sec1]{Introduction},
Theorems \ref{theorem 5.19.1} and \ref{theorem 6.4.1}
are known and even in much stronger forms
for quite some time and we give them with proofs
just to show that there is a probabilistic technique
to derive them and also to prepare some necessary
tools for proving our next result, which
is about Lipschitz continuity of $v$
with constant \emph{independent} of $\delta_{0}$.
As usual, in this case we need the following:
%
%
%as2.2 #&#
\begin{assumption}
\label{assumption 5.29.02}
There exists a $ \delta_{1} \in(0,1]$ such that
for any
$\alpha\in A$, $\beta\in B$, $x\in\bR^{d}$ and
$p\in\bR^{k}$ we have
\[
c^{\alpha\beta} (p,x)\geq\delta_{1}.
\]
\end{assumption}

%
%re2.1 #&#
\begin{remark}
\label{remark 5.29.03}
Assume that $D$ lies in the ball of radius $R$
centered at the origin. For $\mu>0$
define $\Psi(x)=\cosh(\mu R)-\cosh(\mu|x|)+2$. It is
easy to check that for $\mu$
large enough depending only on $\delta_{0},K_{0}$ and $d$,
the function $\Psi$
is infinitely differentiable
on $\bR^{d}$, $\Psi\geq2$ on $D$
and $(L^{\alpha\beta}+c^{\alpha\beta})
\Psi\leq-1$ on $D$ for all $\alpha,\beta$.
This is a so-called global barrier for $D$.

We modify it for $|x|\geq R$ in such a way that
it will be still infinitely differentiable
on $\bR^{d}$, have bounded derivatives and be such that
$\Psi\geq1$ on $\bR^{d}$. We keep the same notation
for the modified function.
By Remark~2.3 of \cite{Kr_2} if we construct
$\check v$ from
\begin{eqnarray*}
\check{\sigma}^{\alpha\beta}(x)&=&\Psi^{1/2}(x) \sigma^{\alpha\beta}(x),\qquad
\check{b}^{\alpha\beta}(x) =\Psi(x)b^{\alpha\beta}(x)+2a^{\alpha\beta}(x)D
\Psi(x),
\\
\check{c}^{\alpha\beta}(x)&=&-L^{\alpha\beta}\Psi(x),\qquad \check{f}^{\alpha\beta}(x)=
f^{\alpha\beta}(x),\qquad \check{g}(x)=\Psi^{-1}(x)g(x),
\end{eqnarray*}
where $D\Psi$ is the gradient of $\Psi$ (a column vector),
in the same way as $v$ was constructed from the original
$\sigma,b,c,f$ and $g$, then $\check v=\Psi^{-1}v$.
By no means the above transformation is
something new; see, for instance,
Sections~1.2 and 2.5 in~\cite{PW}.
Just in case, observe that
now $c^{\alpha\beta}$ influences $\check v$ through
$\check{c}^{\alpha\beta}$, which is bigger than one
(remember that $c^{\alpha\beta}\geq0$).
This shows that
without restricting generality we could have supposed that
Assumption~\ref{assumption 5.29.02} is satisfied
even in Theorems \ref{theorem 5.19.1}
and~\ref{theorem 6.4.1}.
\end{remark}

Introduce
\begin{eqnarray*}
\hat\sigma^{\alpha\beta}(x,y)&=&r^{\alpha\beta}(x,y) \sigma^{\alpha\beta}
\bigl(p^{\alpha\beta}(x,y),x \bigr)P^{\alpha\beta}(x,y),
\\
(\hat a,\hat b,\hat c,\hat f)^{\alpha\beta}(x,y)& =& \bigl[r^{\alpha\beta}(x,y)
\bigr]^{2}(a,b,c,f)^{\alpha\beta} \bigl( p^{\alpha\beta}(x,y),x \bigr),
\end{eqnarray*}
and for
unit $\xi\in\bR^{d}$ introduce a convex
function $\|\sigma\|^{2}_{\xi}$ on the set of
$d\times d_{1}$ matrices by
%
%
%e2.4 #&#
%e2.4 ###
\begin{equation}
\label{5.22.1} \|\sigma\|^{2}_{\xi}:= \|\sigma
\|^{2}-\bigl|\xi^{*}\sigma\bigr|^{2}=\bigl\| \bigl(I-\xi
\xi^{*} \bigr) \sigma\bigr\|^{2}, \qquad \|\sigma\|^{2}=\sum
_{i,j}\sigma^{2}_{ij},
\end{equation}
where $I$ is the unit $d\times d$ matrix.

%
%as2.3 #&#
\begin{assumption}
\label{assumption 5.19.2}
For all $\alpha\in A$, $\beta\in B$
and $x,y\in\bR^{d}$
\[
\delta_{1}^{-1}\geq r^{\alpha\beta}(x,y)\geq
\delta_{1}.
\]
\end{assumption}
%
%
%as2.4 #&#
\begin{assumption}
\label{assumption 5.29.2}
There exist constants $\delta\geq2\delta_{1}$,
$ \varepsilon_{0}>0$ and $\mu\geq1$ such that
for all $\alpha\in A$, $\beta\in B$ and $x,y\in\bR^{d}$,
for which $|x-y|\leq\varepsilon_{0}$, we have
%
%e2.5 ###
\begin{eqnarray}\label{5.24.2}
&&\bigl\|\hat\sigma^{\alpha\beta} (x, y)- \sigma^{\alpha\beta}(y)\bigr\|^{2}_{\xi}
+2 \bigl\langle x-y,\hat b^{\alpha\beta} (x, y) - b^{\alpha\beta}(y) \bigr\rangle
\nonumber
\\[-8pt]
\\[-8pt]
\nonumber
&&\qquad\leq2 \bigl( c^{\alpha\beta}(y) -\delta \bigr)|x-y|^{2}+4
\mu \bigl\langle x-y, a^{\alpha\beta}(x) (x-y) \bigr\rangle,
\end{eqnarray}
where $\xi=(x-y)/|x-y|$.
\end{assumption}

%
%re2.2 #&#
\begin{remark}
\label{remark 5.21.1}
If $d=1$, then
for any $d\times d_{1}$-matrix $\sigma$ and unit $\xi\in\bR^{d}$,
we have $\|\sigma\|=|\xi^{*}\sigma|$, so that
in that case the term involving $\sigma$ in \eqref{5.24.2}
disappears.
Also notice that if $\sigma$ and $b$ are independent of $p$,
and $r\equiv1$, $p\equiv0$, and $P\equiv I$, then
$(\hat a, \hat\sigma,\hat b,\hat c)^{\alpha\beta}(x,y)
=(a,\sigma,b,c)^{\alpha\beta}(x)$, and condition \eqref{5.24.2}
becomes
%
%e2.6 ###
\begin{eqnarray}\label{5.29.5}
&&\bigl\|\sigma^{\alpha\beta} (x )-\sigma^{\alpha\beta}(y)\bigr\|^{2}_{\xi}
+2 \bigl\langle x-y,b^{\alpha\beta} (x ) -b^{\alpha\beta}(y) \bigr\rangle
\nonumber
\\[-8pt]
\\[-8pt]
\nonumber
&&\qquad \leq2 \bigl(c^{\alpha\beta}(y) -\delta \bigr)|x-y|^{2}+4
\mu \bigl\langle x-y, a^{\alpha\beta}(x) (x-y) \bigr\rangle,
\end{eqnarray}
which is satisfied with any $\delta$ on the
account of choosing a sufficiently large $\mu$
(depending on $\delta_{0}$ and $K_{1}$)
since $\sigma$ and $b$ are Lipschitz continuous.
Therefore, Theorem~\ref{theorem 5.19.1}
is a particular case of Theorem~\ref{theorem 5.22.1}.
It is also worth noting that if $d=1$, condition
\eqref{5.29.5} is satisfied with $\mu=0$
when $b^{\alpha\beta}(x)$
are decreasing functions of $x$ and $c^{\alpha\beta}\geq\delta$.

In Section~\ref{section 5.29.6} we give more examples when one can check
Assumption~\ref{assumption 5.19.2}.
\end{remark}

Introduce
\begin{eqnarray*}
&&H \bigl(p,x, u,(u_{i}),(u_{ij}) \bigr)
\\
&&\qquad=\supinf_{\alpha\in A\,  \beta\in B} \bigl[ a^{\alpha\beta}_{ij}(p,x
)u_{ij}+ b_{i}^{\alpha\beta}(p,x )u_{i}-
c^{\alpha\beta}(p,x )u + f^{\alpha\beta}(p,x ) \bigr].
\end{eqnarray*}

%
%as2.5 #&#
\begin{assumption}
\label{assumption 5.22.1}
The set of $(x,u,(u_{i}),(u_{ij}))$ such that
%
%
%e2.7 #&#
%e2.7 ###
\begin{equation}
\label{4.24.3} H \bigl(p,x, u,(u_{i}),(u_{ij}) \bigr)
\leq0
\end{equation}
is \emph{independent} of $p$ and the same is true
if we reverse the sign of the inequality.
\end{assumption}

Note that the next result does not cover
Theorem~\ref{theorem 6.4.1} and by ``the above assumptions''
we mean all assumptions which are stated above in this section.

%
%th2.3 #&#
\begin{theorem}
\label{theorem 5.22.1}
Under the above assumptions also suppose that
either
$D=\bR^{d}$, or
$D$ is bounded and
there exists a nonnegative function $G\in C^{0,1}(\bar{D})
\cap C^{2}_{\mathrm{loc}}(D)$ such that $G=0$ on
$\partial D$ and
\[
L^{\alpha\beta}G(p,x)\leq-1
\]
in $D$ for any $p$.

Then $v$ is Lipschitz continuous in $\bR^{d}$ with Lipschitz
constant \emph{independent} of $\delta_{0}$.
\end{theorem}

%
%re2.3 #&#
\begin{remark}
\label{remark 5.29.1}
If $D$ is bounded and satisfies
the uniform exterior ball condition, the function $G$
always exists since the operators $L^{\alpha\beta}$
are uniformly nondegenerate, have bounded coefficients and
$c^{\alpha\beta}\geq0$. However, the proof of this
well-known fact relies on the uniform nondegeneracy
and gives a function $G$ depending on $\delta_{0}$.
The reader should understand that there are plenty of cases
when this assumption is satisfied,
even for degenerate operators; see, for instance,
Example~\ref{example 5.30.1} with $\delta_{0}=0$.
\end{remark}

Finally, we state one more result, which
was actually the main motivation
of writing the whole series consisting of \cite{Kr_1,Kr_2}
and the present article,
as we have pointed out in the \hyperref[sec1]{Introduction}.
We take $D=\bR^{d}$ and suppose that
all above assumptions are satisfied
and $\sigma,b,c,f$ are independent of $p$.

Set
\[
A_{1}=A,
\]
and let $A_{2}$ be a
separable metric space having no common points with $A_{1}$.

%
%as2.6 #&#
\begin{assumption}
\label{assumption 4.29.2}
The functions $
\sigma^{\alpha\beta}( x)$,
$ b^{\alpha\beta}( x)$, $
c^{\alpha\beta}( x)$ and
$ f^{\alpha\beta}( x)$
are also defined on
$A_{2}\times B \times\bR^{d}$ in such a way that they are
\emph{independent} of $\beta$ and
satisfy Assumptions \ref{assumption 5.19.1}(i), (iii), (iv)
with the same constants $K_{0}$, $K_{1}$
and, of course, with $A_{2}$ in place of $A$.
\end{assumption}

Define
\[
\hat{A}=A_{1}\cup A_{2}.
\]

Then we introduce $\hat{\frA}$ as the set of progressively measurable
$\hat{A}$-valued processes and $\hat{\mathbb{B}}$ as the set
of $\frB$-valued functions $ \bbeta(\alpha_{\cdot})$
on $\hat{\frA}$ such that,
for any $T\in[0,\infty)$ and any $\alpha^{1}_{\cdot},
\alpha^{2}_{\cdot}\in\hat{\frA}$ satisfying
\[
P \bigl( \alpha^{1}_{t}=\alpha^{2}_{t}
\mbox{ for almost all } t\leq T \bigr)=1,
\]
we have
\[
P \bigl( \bbeta_{t} \bigl(\alpha^{1}_{\cdot}
\bigr)= \bbeta_{t} \bigl(\alpha^{2}_{\cdot} \bigr)
\mbox{ for almost all } t\leq T \bigr)=1.
\]

For a constant $K\geq0$, set
\[
v_{K}(x)=\infsup_{\bbeta\in\hat{\mathbb{B}}\,  \alpha_{\cdot}\in
\hat{\frA}} v^{\alpha_{\cdot}\bbeta(\alpha_{\cdot})}_{K}(x),
\]
where
\begin{eqnarray*}
v^{\alpha_{\cdot}\beta_{\cdot} }_{K}(x)&=& E_{x}^{\alpha_{\cdot}\beta_{\cdot} } \int
_{0}^{\infty} f_{K} ( x_{t})e^{-\phi_{t} }
\,dt =:v^{\alpha_{\cdot}\beta_{\cdot} } (x)-K E_{x}^{\alpha_{\cdot}\beta_{\cdot} }\int
_{0}^{\gamma} I_{\alpha_{t}\in A_{2}}e^{-\phi_{t} } \,dt,
\\
f^{\alpha\beta}_{K}( x)&=&f^{\alpha\beta}( x)-KI_{\alpha\in A_{2}}.
\end{eqnarray*}
The above formula extends $v^{\alpha_{\cdot}\beta_{\cdot} } (x)$,
initially defined for $\alpha_{\cdot}\in\frA$ and $\beta_{\cdot}
\in\frB$, on the set $\hat\frA\times\frB$. Of course,
\eqref{2.12.2} is preserved with $\tau=\infty$, and no $g$ is
involved.

%
%th2.4 #&#
\begin{theorem}
\label{theorem 6.3.1}
There is a constant $N$,
depending only on the constants in all above assumptions
(but not on $K$),
such that $|v_{K}(x)-v(x)|\leq
N/K$ for all $x\in\bR^{d}$ and $K\geq1$.
\end{theorem}

%
%re2.4 #&#
\begin{remark}
\label{remark 6.3.1}
In one of the main cases of interest $v_{K}$ turns out to have
second-order derivatives bounded by a constant times $K$
if $K\geq1$; see Section~7 in~\cite{Kr_2}. From the point of
view of finite-difference approximations it should be easier
to approximate ``smooth'' functions $v_{K}$ than $v$. However, the author
has no idea how to prove a fact similar to Theorem~\ref{theorem 6.3.1} for finite-difference equations.

In this connection it would be very interesting to find any
proof of Theorem~\ref{theorem 6.3.1} not using probability theory,
of course, defining $v_{K}$ and $v$ as viscosity solutions
of the corresponding Isaacs equations.
\end{remark}

%s3 #&#
%s3 ###
\section{Comments and examples}
\label{section 5.29.6}
%
%
%re3.1 #&#
\begin{remark}
\label{remark 5.29.2}
Let $\sigma$ and $b$ be independent of $\alpha$ and $\beta$,
and consider a particular case where $d_{1}=d$, and
equation \eqref{5.11.1} is
%
%
%e3.1 #&#
%e3.1 ###
\begin{equation}
\label{3.8.4} x_{t}=x+\int_{0}^{t}
\sigma(x_{s}) \,dw_{s},
\end{equation}
where $\sigma$ is an $\cO$-valued Lipschitz continuous function.
Then the left-hand side of~\eqref{5.24.2} vanishes
for $r\equiv1$ and $P(x,y)=\sigma^{*}(x)\sigma(y)$. Of course,
this is not a big surprise since $x_{t}$ is just a Brownian
motion starting at $x$. Still one can see that the parameters $P$
take care of rotations of the
increments of the original Wiener process and basically show that
\eqref{5.24.2} is a condition on $a$ rather than $\sigma$.
\end{remark}

%
%re3.2 #&#
\begin{remark}
\label{remark 6.3.4}
The function $v$ will not change if we change
$\sigma,b,c,f$ outside~$D$. In connection with this
it is worth noting that in Assumption~\ref{assumption 5.29.2} we may restrict $x$ and $y$ to
$D_{\varepsilon_{0}}$ which is the $\varepsilon_{0}$
neighborhood of $D$. Indeed, if only thus restricted
Assumption~\ref{assumption 5.29.2} is satisfied
we could just change $c$ outside $D$ so that it will be bigger than
the original one and become any large constant
outside $D_{\varepsilon_{0}}$. Then Assumption~\ref{assumption 5.29.2}
will be satisfied in the form it is stated.
\end{remark}

%
%re3.3 #&#
\begin{remark}
\label{remark 5.21.01}
For later discussion
we show that Assumption~\ref{assumption 5.29.2} can be replaced
with a slightly more transparent one. We will be only
concerned with Assumption~\ref{assumption 5.29.2}
leaving other assumptions aside.

Denote by $\cSk$ the set of
$d_{1}\times d_{1}$ skew-symmetric
matrices and
assume that
for each $\alpha\in A$, $\beta\in B$
and $\xi\in\bR^{d}$, the following
functions on $\bR^{d}$ are also given:
$\cSk$-valued
function $\Theta^{\alpha\beta}(x,\xi)$, $k\times d$ matrix-valued
function $p^{\alpha\beta}(x)$, and $\bR^{d}$-valued
function $r^{\alpha\beta}(x)$.

For a differentiable function $u(p,x)$ and $\xi\in\bR^{d}$, introduce
\[
\partial_{\xi}u^{\alpha\beta}(x )=\xi_{i}u_{x_{i}}
(0,x) + \bigl(p^{\alpha\beta}(x)\xi \bigr)_{j}u_{p_{j}}
(0,x).
\]
Also denote $\operatorname{Conv} (D)$ the open convex hull of $D$.
%
%
%as3.1 #&#
\begin{assumption}
\label{assumption 5.30.1}
(i) For $|\xi|\leq1$ the above functions are bounded by
$K_{0}$ and
$\Theta^{\alpha\beta}(x,y)$ is a linear function of $y$
[in particular $\Theta^{\alpha\beta}(x,0)=0$].\vspace*{-6pt}
\begin{longlist}[(iii)]
\item[(ii)] For any $\alpha\in A$ and $\beta\in B$ the functions
$\sigma^{\alpha\beta}(p,x)$ and $b^{\alpha\beta}(p,x)$
are continuously differentiable with respect to $(p,x)
\in\bR^{k}\times\bR^{d}$,
and their first-order derivatives are bounded by
$K_{1}$. Furthermore, their derivatives are uniformly
continuous with respect to $(p,x)$ uniformly with respect
to $(\alpha,\beta)\in A\times B $.

\item[(iii)] There are constants $\mu\geq1$ and $\delta\geq2\delta_{1}$
such that for any unit $\xi\in\bR^{d}$ and
$(\alpha,\beta,x)\in A\times B\times\operatorname{Conv} (D)$, we have
%
%e3.2 ###
\begin{eqnarray}\label{5.21.1}
&&\bigl\| \partial_{\xi} \sigma^{\alpha\beta} (x) + \bigl\langle
r^{\alpha\beta}(x),\xi \bigr\rangle \sigma^{\alpha\beta} (x)+\sigma^{\alpha\beta}
(x) \Theta^{\alpha\beta}(x,\xi) \bigr\|^2_{\xi}\nonumber
\\
&&\quad{} + 2 \bigl\langle\xi,\partial_{\xi}b^{\alpha\beta}(x)+ 2
\bigl\langle r^{\alpha\beta}(x),\xi \bigr\rangle b^{\alpha\beta}(x) \bigr\rangle
\\
&&\qquad\le2 \bigl( c^{\alpha\beta}(x) - \delta_{1}- \delta \bigr)+ 4\mu \bigl
\langle\xi,a^{\alpha\beta}(x)\xi \bigr\rangle.\nonumber
\end{eqnarray}
\end{longlist}
\end{assumption}
Introduce
%
%e3.3 ###
\begin{eqnarray}\label{5.30.2}
r^{\alpha\beta}(x,y)&=&1+ \bigl\langle r^{\alpha\beta}(y),x-y \bigr\rangle,\qquad
p^{\alpha\beta}(x,y)=p^{\alpha\beta}(y) (x-y),
\nonumber
\\[-8pt]
\\[-8pt]
\nonumber
 P^{\alpha\beta}(x,y) &=&\exp\Theta^{\alpha\beta}(y,x-y).
\end{eqnarray}

We claim that there exists an $\varepsilon_{0}>0$,
depending only on $K_{0},K_{1},\delta_{1},d$,
and the moduli of continuity in $(p,x)$ of
the derivatives of
$\sigma^{\alpha\beta}(p,x)$ and $b^{\alpha\beta}(p,x)$
with respect to
$(p,x)$, such that Assumption~\ref{assumption 5.29.2} is satisfied
with $x,y$ restricted to $D$.

To prove the claim, fix $y\in D$ and a unit $\xi\in\bR^{d}$,
and for $t\geq0$ introduce $x(t)=y+t\xi$, so that \eqref{5.24.2}
becomes
%
%e3.4 ###
\begin{eqnarray}\label{5.22.2}
&&\bigl\|\hat\sigma^{\alpha\beta} \bigl(x(t), y \bigr)-\sigma^{\alpha\beta}(y)
\bigr\|^{2}_{\xi} + 2t \bigl\langle\xi,\hat b^{\alpha\beta}
\bigl(x(t), y \bigr) -b^{\alpha\beta}(y) \bigr\rangle
\nonumber
\\[-8pt]
\\[-8pt]
\nonumber
&&\qquad \leq2 \bigl(c^{\alpha\beta}(y) -\delta \bigr)t^{2}+4\mu
\bigl\langle\xi, a^{\alpha\beta} \bigl(x(t) \bigr)\xi \bigr\rangle
t^{2},
\end{eqnarray}
which we want to prove for $t\in(0,\varepsilon_{0}]$.
For simplicity of notation we will drop the superscripts
$\alpha,\beta$ in a few lines below.

Observe that
\[
\hat\sigma \bigl(x(t), y \bigr)-\sigma(y)=\int_{0}^{t}
\xi_{i} \hat\sigma_{x_{i}} \bigl(x(s), y \bigr) \,ds,
\]
where
\begin{eqnarray*}
&&\xi_{i}\hat\sigma_{x_{i}} \bigl(x(s), y \bigr)\\
&&\qquad=\bigl\langle
r (y),\xi \bigr\rangle \sigma \bigl(sp (y)\xi,x(s) \bigr) P \bigl(x(s),y \bigr)
\\
&&\qquad\quad{}+r \bigl(x(s),y \bigr) \bigl[\xi_{i}\sigma_{x_{i}} \bigl(sp
(y)\xi,x(s) \bigr) + \bigl(p (y)\xi \bigr)_{j} \sigma_{p_{j}}
\bigl(sp (y)\xi,x(s) \bigr) \bigr]\\
&&\qquad\quad{}\times P \bigl(x(s),y \bigr)
\\
&&\qquad\quad{}+r \bigl(x(s),y \bigr) \sigma \bigl(sp (y)\xi,x(s) \bigr) \Theta(y,\xi)P
\bigl(x(s),y \bigr)
\\
&&\qquad=: \bigl\langle r (y),\xi \bigr\rangle \sigma(y) +\partial_{\xi}
\sigma(y) +\sigma(y) \Theta(y,\xi)+R (s ),
\end{eqnarray*}
and $R (s )$ is introduced by the above
equality.

Owing to the convexity of function \eqref{5.22.1}
and Assumption~\ref{assumption 5.30.1},
there exists an $\varepsilon_{0}>0$ such that
for all $t\in(0,\varepsilon_{0}]$ and all values of other arguments,
we have
\begin{eqnarray*}
&&\bigl\|\hat\sigma^{\alpha\beta} \bigl(x(t), y \bigr)-\sigma^{\alpha\beta}(y)
\bigr\|^{2}_{\xi}- 4\mu \bigl\langle\xi, a^{\alpha\beta} \bigl(x(t)
\bigr)\xi \bigr\rangle t^{2}
\\
&&\qquad\leq t^{2}\bigl \|\partial_{\xi} \sigma^{\alpha\beta}(y) + \bigl
\langle r^{\alpha\beta}(y),\xi \bigr\rangle \sigma^{\alpha\beta}(y) +
\sigma^{\alpha\beta}(y) \Theta^{\alpha\beta}(y,\xi)\bigr\|^{2}_{\xi}
\\
&&\qquad\quad{}-4\mu \bigl\langle\xi, a^{\alpha\beta}(y)\xi \bigr\rangle t^{2}+t^{2}
\delta_{1}.
\end{eqnarray*}

It is even easier to prove that, by reducing $\varepsilon_{0}$
if necessary, we have that for
$t\in(0,\varepsilon_{0}]$ and all values of other arguments
\begin{eqnarray*}
&&t \bigl\langle\xi,\hat b^{\alpha\beta} \bigl(x(t), y \bigr) -b^{\alpha\beta}(y)
\bigr\rangle
\\
&&\qquad\leq t^{2} \bigl\langle\xi, \partial_{\xi}
b^{\alpha\beta}(y)+ 2 \bigl\langle r^{\alpha\beta}(y),\xi \bigr\rangle
b^{\alpha\beta}(y) \bigr\rangle+t^{2}\delta_{1}.
\end{eqnarray*}

Hence, by assumption, the left-hand side of
\eqref{5.22.2} is less than
\[
t^{2} \bigl[2 \bigl( c^{\alpha\beta}(y) -\delta_{1}-\delta
\bigr)+ 4\mu \bigl\langle\xi,a^{\alpha\beta}(y)\xi \bigr\rangle \bigr] +
2t^{2}\delta_{1},
\]
which is the right-hand side of \eqref{5.22.2}.
\end{remark}
%
%
%re3.4 #&#
\begin{remark}
\label{remark 5.21.2}
Consider the case that $\sigma$ and $b$ are
independent of $\alpha$ and $\beta$.
Let $d =1$, $c>0$ and $D=(-1,1)$.
Assume that $a=a_{0}+\delta_{0}$,
where $a_{0}\geq0$. In that case,
as it follows from the
arguments in Remarks \ref{remark 5.21.1}
and \ref{remark 5.21.01},
we do not need to assume
that $\sigma'$ is continuous.
We still assume that $a$, $b'$ and $c$ are continuous.
Then by Remark~\ref{remark 5.21.1}
Assumption~\ref{assumption 5.29.2}
is satisfied with $\mu$ depending on $\delta_{0}$, among other
things.

However,
assume additionally that at every point $x\in[-2,2]$ where
\[
a_{0}(x)=b(x)=0
\]
we have
%
%
%e3.5 #&#
%e3.5 ###
\begin{equation}
\label{5.29.04} b'(x) < c(x).
\end{equation}

We claim that then Assumption~\ref{assumption 5.29.2}
is satisfied with $x,y$ restricted
to $[-2,2]$ with some $\delta,\delta_{1},\varepsilon_{0}$,
and $\mu$ \emph{independent} of $\delta_{0}$ and hence,
by Remark~\ref{remark 6.3.4}, it will be satisfied in the original form,
making the assertion of Theorem~\ref{theorem 5.22.1} valid
in case $D=(-1,1)$.

To prove the claim, we use Remark~\ref{remark 5.21.01} and
observe that for $r=-nb/2$,
$\delta_{1}+\delta=1/n$, $\mu= n$ and $|\xi|=1$ condition \eqref{5.21.1}
is satisfied if
%
%
%e3.6 #&#
%e3.6 ###
\begin{equation}
\label{5.29.05} b'(x) \le c(x) - \frac{1}n + n \bigl(
a_{0}(x) + \bigl| b(x)\bigr |^2 \bigr).
\end{equation}

Suppose that for any $n=1,2,\ldots$ we can find
a point $x_n\in[-2,2]$ at which the inequality converse to
\eqref{5.29.05} holds. Then we can extract from the sequence
$x_{n}$ a subsequence that
converges to an
$ x_0 \in[-2,2] $. Clearly, for large $n$,
\[
a_{0}(x_n) + \bigl| b(x_n) \bigr|^2 \le
Nn^{-1},
\]
where $N=\sup b' +1$. Therefore,
$a_{0}(x_0) +  | b(x_0)  |^2 =0$ and
\[
b'(x_n) \ge c(x_n) - 1/n,\qquad
b'(x_0) \ge c(x_0).
\]
We have obtained a contradiction to \eqref{5.29.04}, so
inequality \eqref{5.29.05} holds in $[-2,2]$ for some $n$
independent of $\delta_{0}$ thus proving our claim.
\end{remark}

%
%ex3.1 #&#
\begin{example}
\label{example 5.30.1}
Consider the one-dimensional equation
%
%
%e3.7 #&#
%e3.7 ###
\begin{equation}
\label{5.30.4} \delta_{0} v''+bx
v'-v=0
\end{equation}
on $[-1,1]$ with data $1$ at $\pm1$,
where constant $b>0$.
This is, of course, a simple example of the Isaacs
equation in a differential ``game''
with the value function $v$. Here the assumption
stated in Theorem~\ref{theorem 5.22.1} concerning $G$ is satisfied
with $G(x)=(1-x^{2})\max(1,1/(2b))$.

If we assume that the solution $v=v_{\delta_{0}}$
admits an estimate of its Lipschitz constant independent
of $\delta_{0}$, then, as is easy to understand, say from
the probabilistic representation of $v_{\delta}$,
the function
\[
v_{0}(x)=Ee^{-\tau_{x}}
\]
would be Lipschitz continuous, where $\tau_{x}$ is the first
exit time of the solution of
\[
x_{t}=x+ \int_{0}^{t}b
x_{s} \,ds
\]
from $(-1,1)$. Since $x_{t}=xe^{bt}$, $\tau_{x}=-b^{-1}\ln|x|$
for $|x|<1$
and $v_{0}(x)=|x|^{1/b}$, which is Lipschitz continuous
only if $b\leq1$.

This example shows that in the situation of Remark~\ref{remark 5.21.2},
if one has $b'(x)>c(x)$ at least at
one point at which $a_{0}(x)
=b(x)=0$, the assertion of Theorem~\ref{theorem 5.22.1} may be no longer
true. In this respect, requiring condition \eqref{5.29.04} at those points
is close to being optimal and it is, actually, necessary for $v$
to be \emph{continuously} differentiable.
\end{example}

%
%re3.5 #&#
\begin{remark}
\label{remark 5.31.4}
Barles in \cite{Ba91} derived first-order derivatives estimates
for viscosity solutions of nonlinear equations
\[
H \bigl( x, u,Du,D^{2}u \bigr)=0
\]
in domains, where $Du=(D_{i}u)$ is the gradient of $u$, and $D^{2}u=
(D^{2}_{ij}u)$
is its Hessian. Our value functions are viscosity solutions
of the corresponding Isaacs equations.
This fact
is proved
in \cite{Ko09} for bounded domains.
The Isaacs equations in this paper are included
in the framework of \cite{Ba91} and many of the equations
in \cite{Ba91} do not fit into our scheme.
Yet it is worth comparing our conditions with the ones from~\cite{Ba91} in the \emph{simplest} example
of \emph{linear} equations with
\[
H \bigl(x,u_{0},u',u'' \bigr)
=a_{ij}(x)u''_{ij}+b_{i}(x)u'_{i}-c(x)u
_{0}+f(x)
\]
for which solutions have probabilistic representations
(with no $\alpha$ and $\beta$ involved).

One of the assumptions in \cite{Ba91}
reads as follows: For any $R>0$ and all large enough $L$,
%
%e3.8 ###
\begin{eqnarray}\label{2.8.1}
&&c\sum_{i=1}^{d}\bigl|u'_{i}\bigr|^{2}
+g\tr u''au''
\nonumber\\
&&\hspace*{20pt}{}- \bigl[u'_{k}D_{k}a_{ij}u''_{ij}+u'_{k}D_{k}b_{i}(x)u'_{i}
-u'_{k}D_{k}c(x)u _{0}
+u'_{k}D_{k}f(x) \bigr]
\\
&&\quad\qquad\geq h,\nonumber
\end{eqnarray}
where $g,h>0$ are some constants $>0$, provided that
%
%
%e3.9 #&#
%e3.9 ###
\begin{equation}\qquad\hspace*{4pt}
\label{2.8.2} |u_{0}|\leq R,\qquad \sum_{i=1}^{d}\bigl|u'_{i}\bigr|^{2}
\geq L,\qquad H \bigl(x,u_{0},u',u''
\bigr)=0, \qquad u''_{ij}u'_{j}=0\ \forall i.\hspace*{-3pt}
\end{equation}
If $c\equiv0$, $b\equiv0$, and both
$f $ and $Df $ vanish at a point $x_{0}$, so that $H(x_{0},0)=0$, then
for $u''=0$ inequality
\eqref{2.8.1} at $x_{0}$ becomes $0\geq h$, which cannot hold
even in the one-dimensional case. Therefore, the one-dimensional
equation
\[
D^{2}u+x^{2}=0
\]
in $(-1,1)$ with zero boundary condition
does not fit in the scheme of \cite{Ba91}.

Equation
\[
\delta_{0}D^{2}u+ (b_{1} x+b_{0})
Du-cu+x^{2}=0
\]
in $(-1,1)$ with zero boundary condition
and constant $c>0,b_{0},b_{1}$ does not fit in
either if $c\leq b_{1}$.

Indeed, if we take $x=0$, $u''=0$, $u_{0} =0$, and $u'$ bigger
by magnitude than $L$, \eqref{2.8.1} becomes
\[
(c-b_{1})\bigl|u'\bigr|^{2}\geq h,
\]
which for large $|u'|$ can only hold if $ b_{1}<c$.
Remark~\ref{remark 5.21.2} shows that one always has an estimate
of the Lipschitz constant of $v$. This estimate is
even \emph{independent}
of $\delta_{0}$, provided that either $
b_{1}x+b_{0}\ne0$ for $x\in[-1,1]$
or $ b_{1}<c$.

It looks like
the methods of \cite{Ba91}
are not adapted to use uniform nondegeneracy and even in the above
examples
lead to the requirement that $c$ be sufficiently large.
\end{remark}

%
%re3.6 #&#
\begin{remark}
\label{remark 6.1.1}
Above we saw that the parameters $\mu$, $r$ and $P$ can play a role
while checking Assumption~\ref{assumption 5.29.2}. We now show how
the external parameters $p$ can be used. Here we consider the situation
in which $\sigma$, $b$, $c$ and $f$
depend only on $x$ and $\alpha$
so that we are dealing with controlled diffusion processes
rather than differential games. Our interest is in obtaining
estimates independent of $\delta_{0}$, and therefore,
from the start in this remark
we focus on \emph{degenerate} processes.

Let $A=\bR$ and consider a
one-dimensional process defined by the equation
%
%
%e3.10 #&#
%e3.10 ###
\begin{equation}
\label{6.1.3} x_{t}=x+\int_{0}^{t}
\sigma(x_{s}) \,dw_{s}+ \int_{0}^{t}
\tanh(x_{s}+2\cos\alpha_{s}) \,ds,
\end{equation}
where $w_{t}$ is a one-dimensional Wiener process,
$\sigma(x)$ is a smooth nonnegative even function satisfying $\sigma(x)
>0$ for $x\in(1,3)$ and vanishing outside $(1,3)$
(and $\alpha_{t}$ is a progressively measurable
$A$-valued process). We also take a sufficiently regular
function $c(x)\geq\delta_{2}$ (independent of $\alpha$
and $\beta$),
where $\delta_{2}>0$, and take $D=\bR$.

If we want to satisfy \eqref{5.21.1}
for $|x|\notin[1,3]$ with $r(x)=0$ (and $\Theta\equiv0$
for having no other options)
and some $\delta$'s we obviously need to have
%
%
%e3.11 #&#
%e3.11 ###
\begin{equation}
\label{6.1.2}\qquad c(x)>1\qquad \mbox{for } |x|\leq1, \qquad c(x)> \cosh^{-2}\bigl(|x|-2\bigr)
\qquad\mbox{for } |x|\geq3.
\end{equation}
The inequalities in \eqref{6.1.2} extend for
$|x|\notin(1+\varepsilon,3-\varepsilon)$ with some $\varepsilon
>0$, and one can find $\mu\geq1$ such that
\eqref{5.21.1} is satisfied
(with some
$\delta$'s) for $|x|\in(1+\varepsilon,3-\varepsilon)$
with $r(x)=0$. Therefore, if we do not use parameter $r$,
then \eqref{5.21.1} reduces to \eqref{6.1.2}.

However, if we take
%
%
%e3.12 #&#
%e3.12 ###
\begin{equation}
\label{6.1.4} r^{\alpha}(x)=-2I_{|x+2\cos\alpha|>\varepsilon} \sinh^{-1}(2x+4
\cos\alpha),
\end{equation}
then the left-hand side of \eqref{5.21.1} becomes
\[
2I_{|x+2\cos\alpha|\leq\varepsilon}\cosh^{-2} (x+2\cos\alpha)\leq2I_{|x+2\cos\alpha|\leq\varepsilon},
\]
and for $|x|\notin(1+\varepsilon,3-
\varepsilon)$ this is strictly less than $2c(x)$ if
%
%
%e3.13 #&#
%e3.13 ###
\begin{equation}
\label{6.1.5} c(x)>1\qquad \mbox{for } |x|\leq1+\varepsilon.
\end{equation}
Hence, with the so specified $r^{\alpha}$
condition, \eqref{5.21.1} reduces to \eqref{6.1.5},
which is a significant improvement over \eqref{6.1.2}.

Next we take $f$ independent of $\alpha$, say $f\equiv1$,
and instead of
\[
b^{\alpha}(x)=\tanh(x +2\cos\alpha)
\]
consider
\[
b^{\alpha}(p,x)=\tanh \bigl(x +2\cos(\alpha+p) \bigr),
\]
where $p\in\bR$. Obviously, Assumption~\ref{assumption 5.22.1} is satisfied.

Take $r^{\alpha}(x)$ from \eqref{6.1.4} and
%
%
%e3.14 #&#
%e3.14 ###
\begin{equation}
\label{6.1.6} p^{\alpha}(x)= (1/2)I_{|x+2\cos\alpha|\leq\varepsilon}I_{|\sin
\alpha|
>\varepsilon}
\sin^{-1}\alpha.
\end{equation}
Then the left-hand side of \eqref{5.21.1} becomes
\begin{eqnarray*}
&& 2I_{|x+2\cos\alpha|\leq\varepsilon} \cosh^{-2} (x+2\cos\alpha)-2I_{|x+2\cos\alpha|\leq\varepsilon}I_{|\sin
\alpha|
>\varepsilon}
\cosh^{-2} (x+2\cos\alpha)
\\
&&\qquad=2I_{|x+2\cos\alpha|\leq\varepsilon}I_{|\sin\alpha|
\leq\varepsilon}\cosh^{-2} (x+2\cos\alpha)
\leq2I_{|x+2\cos\alpha|\leq\varepsilon} I_{|\sin\alpha|
\leq\varepsilon},
\end{eqnarray*}
and the latter is zero if $|x|\leq1+\varepsilon$
and $\varepsilon$ is sufficiently small. Thus adding $p^{\alpha}(x)$
into the picture eliminates condition \eqref{6.1.5}
entirely, and there is nothing more than $
c(x)\geq\delta_{2}$ required of $c(x)$ in order for \eqref{5.21.1}
to be satisfied
with $r^{\alpha}(x)$ from \eqref{6.1.4} and $p^{\alpha}(x)$
from \eqref{6.1.6}.

By the way, the Isaacs (Bellman) equation in this case is
\[
a(x)D^{2}v(x)+ \bigl(Dv(x) \bigr)\tanh \bigl[x+2\operatorname{sign}
\bigl(Dv(x) \bigr) \bigr]-c(x)v(x)+f(x)=0,
\]
where $a=(1/2)\sigma^{2}$. This equation
suggests a different representation of the value
function with $A=\{\pm1\}$ when using parameters $p$
becomes unnecessary (and impossible) but using $r$ will suffice.
In this connection it is worth mentioning that
much more sophisticated use of the external
parameters $p$ can be found in \cite{Kr89},
where in an example of (degenerate)
complex Monge--Amp\`ere equation
they are shown to be indispensable in proving the
global $C^{1,1}$
regularity of solutions.
\end{remark}

%s4 #&#
%s4 ###
\section{Some underlying ideas}
\label{section 3.8.1}
This article is written for probabilists
and the translation of the proof of the central
Theorem~\ref{theorem 6.3.1}
in PDE terms or in terms of the theory
of viscosity solutions is unknown to the author.
On the other hand, such a translation may exist
for Theorem~\ref{theorem 5.22.1} and the interested,
more PDE oriented, reader can find
in Section~8.5 of
\cite{Kr_4} analytical tools allowing one to prove
an analog of Theorem~\ref{theorem 5.22.1}
for Bellman's equations.

However, for probabilists the following explanation
of ideas behind the proof of Theorem~\ref{theorem 5.22.1}
might be helpful.
The main idea is that while differentiating $v(x)$
with respect to $x$ we can take different
representation for $v$ at different points.
We explain how various terms in \eqref{5.21.1}
appear naturally on two examples
of stochastic equations without games.

%
%ex4.1 #&#
\begin{example}
\label{example 3.9.1}
In Remark~\ref{remark 5.29.2}, take
a smooth bounded $f(x)$ and define
%
%
%e4.1 #&#
%e4.1 ###
\begin{equation}
\label{3.8.3} v(x)=E \int_{0}^{\infty}e^{-t}f
\bigl(x^{x}_{t} \bigr) \,dt,
\end{equation}
where we use the same stipulation about indices
as before and do not write $\alpha$ and $\beta$ because
nothing is depending on these parameters.
One can formally differentiate $v(x)$ and obtain that
for any $\xi\in\bR^{d}$,
%
%
%e4.2 #&#
%e4.2 ###
\begin{equation}
\label{3.8.05} v_{(\xi)}(x)=E \int_{0}^{\infty}e^{-t}f_{(\xi_{t})}
\bigl(x^{x}_{t} \bigr) \,dt,
\end{equation}
where $\xi_{t}$ is defined as the solution of
\[
d\xi_{t}=\sigma_{(\xi_{t})} \bigl(x^{x}_{t}
\bigr) \,dw_{t},\qquad \xi_{0}=\xi.
\]
Actually, it is not hard to see that \eqref{3.8.3}
is indeed true, provided that
%
%
%e4.3 #&#
%e4.3 ###
\begin{equation}
\label{3.8.5} E|\xi_{t}|\leq Ne^{ \gamma t},
\end{equation}
where $N$ is a constant and a constant $\gamma<1$.
In that case the right-hand side of~\eqref{3.8.05} is well defined. This may not
happen if the derivatives of $\sigma$ are big.

However, observe that for any $d\times d$-valued
skew-symmetric progressively
measurable process $\Theta_{t}$ and any $\varepsilon$
we also have
%
%
%e4.4 #&#
%e4.4 ###
\begin{equation}
\label{3.8.04} v(x+\varepsilon\xi)=E \int_{0}^{\infty}e^{-t}f
\bigl(x^{x}_{t}(\varepsilon) \bigr) \,dt,
\end{equation}
where $x^{x}_{t}(\varepsilon)$ is defined as a unique
solution of
\[
dx_{t}=\sigma(x_{t})e^{\varepsilon\Theta_{t}} \,dw_{t},\qquad
x_{0}=x.
\]
Formula \eqref{3.8.04} is indeed true because
\[
e^{\varepsilon\Theta_{t}} \,dw_{t}=db_{t},
\]
where $b_{t}$ is a Wiener process and the distributions
of solutions of \eqref{3.8.4}
are independent of which
Wiener process is involved. Now let us formally differentiate
\eqref{3.8.04} through with respect to $\varepsilon$
at $\varepsilon=0$. We again obtain
\eqref{3.8.05}, but this time $\xi_{t}$ satisfies
%
%
%e4.5 #&#
%e4.5 ###
\begin{equation}
\label{3.8.6} d\xi_{t}= \bigl[\sigma_{(\xi_{t})}
\bigl(x^{x}_{t} \bigr)+\sigma \bigl(x^{x}_{t}
\bigr) \Theta_{t} \bigr] \,dw_{t}, \qquad \xi_{0}= \xi.
\end{equation}
Here the coefficient of $dw_{t}$ vanishes if we take
$\Theta_{t}=-\sigma^{*}(x^{x}_{t})\sigma_{(\xi_{t})}(x^{x}_{t})$,
so that $\xi_{t}\equiv\xi$ and nothing like
\eqref{3.8.5} is an issue any longer. The reader may object that
one cannot take $\Theta_{t}=-
\sigma^{*}(x^{x}_{t})\sigma_{(\xi_{t})}(x^{x}_{t})$
before solving \eqref{3.8.6}. Then take
$\Theta_{t}=-\sigma^{*}(x^{x}_{t})\sigma_{(\xi)}(x^{x}_{t})$
and use that $\xi_{t}\equiv\xi$ satisfies \eqref{3.8.6}.

For any $\Theta_{t}$ we have from \eqref{3.8.6} that
\[
d|\xi_{t}|^{2}=\bigl\|\sigma_{(\xi_{t})}
\bigl(x^{x}_{t} \bigr)+\sigma \bigl(x^{x}_{t}
\bigr) \Theta_{t}\bigr\|^{2} \,dt+dm_{t},
\]
where $m_{t}$ is a local martingale.
This shows the origin of $\sigma^{\alpha\beta} (x)
\Theta^{\alpha\beta}(x,\xi)$ in \eqref{5.21.1}.
The subscript $\xi$ appears there after we compute
$d|\xi_{t}|$.
%%%
\end{example}

%
%ex4.2 #&#
\begin{example}
\label{example 3.9.2}
Consider the one-dimensional It\^o equation
\[
dx_{t}=\sigma(x_{t}) \,dw_{t}+b(x_{t})
\,dt, \qquad x_{0}=x
\]
with one-dimensional $w_{t}$,
and introduce $v(x)$ as in \eqref{3.8.3}, so that $c=1$.
Then we again
have \eqref{3.8.05} provided that \eqref{3.8.5} holds
with a $\gamma<1$ and $\xi_{t}$ defined as a unique
solution of
%
%
%e4.6 #&#
%e4.6 ###
\begin{equation}
\label{3.9.4} d\xi_{t}=\sigma_{(\xi_{t})} \bigl(x^{x}_{t}
\bigr) \,dw_{t} +b_{(\xi_{t})} \bigl(x^{x}_{t}
\bigr) \,dt,\qquad \xi_{0}=\xi.
\end{equation}
The solution of \eqref{3.9.4} is known to be
\[
\xi_{t}=\xi m_{t}\exp\int_{0}^{t}b'
\bigl(x^{x}_{s} \bigr) \,ds,
\]
where
\[
m_{t}=\exp \biggl(\int_{0}^{t}
\sigma' \bigl(x^{x}_{s} \bigr)
\,dw_{s} -(1/2)\int_{0}^{t}\bigl|
\sigma' \bigl(x^{x}_{s} \bigr)\bigr|^{2}
\,ds \biggr)
\]
is at least a supermartingale. Hence \eqref{3.8.5} becomes
\[
Em_{t}\exp\int_{0}^{t}b'
\bigl(x^{x}_{s} \bigr) \,ds\leq Ne^{ \gamma t}
\]
and a sufficient condition for that to happen is $b'\leq
\gamma c$ (since $Em_{t}\leq1$).

However, one can use a random time change and get
a different representation for $v$. Namely, take
any progressively measurable real-valued bounded process
$r_{t}$ and for $\varepsilon$ such that
$1+2\varepsilon r_{t}\geq1/2$ introduce $x^{x}_{t}(\varepsilon)$
as a unique solution of
%
%
%e4.7 #&#
%e4.7 ###
\begin{equation}
\label{3.9.6} dx_{t}=\sqrt{1+2\varepsilon r_{t}}
\sigma(x_{t}) \,dw_{t} +(1+2\varepsilon r_{t})b(x_{t})
\,dt,\qquad x_{0}=x.
\end{equation}
Then it is well known that
%
%
%e4.8 #&#
%e4.8 ###
\begin{equation}
\label{3.9.5} v(x)=E\int_{0}^{\infty}f
\bigl(x^{x}_{t}(\varepsilon) \bigr) (1+2\varepsilon
r_{t}) \exp \biggl(-\int_{0}^{t}(1+2
\varepsilon r_{s}) \,ds \biggr) \,dt.
\end{equation}
We substitute $x+\varepsilon\xi$ in place of $x$ in
\eqref{3.9.5} and differentiate with respect to $\varepsilon$
at $\varepsilon=0$. Then instead of \eqref{3.8.05}
we obtain
%
%
%e4.9 #&#
%e4.9 ###
\begin{equation}
\label{3.9.8} v_{(\xi)}(x)=E\int_{0}^{\infty}
\biggl[f_{(\xi_{t})} \bigl(x^{x}_{t} \bigr)
+2r_{t}f \bigl(x^{x}_{t} \bigr)-2f
\bigl(x^{x}_{t} \bigr)\int_{0}^{t}r_{s}
\,ds \biggr] e^{-t} \,dt,
\end{equation}
where $\xi_{t}$ is defined by the equation
%
%
%e4.10 #&#
%e4.10 ###
\begin{equation}
\label{3.9.7} d\xi_{t}=[\sigma_{(\xi_{t})}+r_{t}
\sigma] \bigl(x^{x}_{t} \bigr) \,dw_{t}+[b_{(\xi_{t})}+2r_{t}b]
\bigl(x^{x}_{t} \bigr) \,dt, \qquad\xi_{0}=\xi.
\end{equation}
After formula \eqref{3.9.8} is obtained for bounded
processes $r_{t}$, it can be extended for a wider class
and we plug $r_{t}=\xi_{t}\alpha(x^{x}_{t})$, where $\alpha(x)$
will be specified later, into~\eqref{3.9.7} solve it
and use the solution in \eqref{3.9.8}. Similarly to what
was said before, these manipulations can be easily
justified if
\[
b'+2\alpha b\leq\gamma c.
\]
This is what \eqref{5.21.1} becomes in our case
with $\mu=0$.

We described the way how the parameters $\Theta$ and $r$ appear.
One can also use a change of probability measure
based on Girsanov's theorem and then one includes
in \eqref{5.21.1} an
additional helping term $(a\xi,\xi)$ with as big
factor as one likes.

More details in a more difficult case of controlled
diffusion processes can be found in \cite{Zh12}.
Note that in the above explanation in both cases in
\eqref{3.8.6} and
\eqref{3.9.7} we first found $\Theta$ and $r$ in the form we like, then
solved these equations and used thus specified $\Theta$
and $r$ in \eqref{3.8.04} and \eqref{3.9.8}.
The same procedure works for controlled diffusion
processes because it is known that one can use any
progressively measurable $\Theta$ and $r$ without
affecting the value function. This property is
unknown, however, for stochastic differential games.
We can only use $\Theta=\Theta(x_{t})$ and $r=r(x_{t})$,
which would not lead to any good result even in the
above examples where $\Theta$ and $r$ depend linearly on $\xi$.
Therefore, what we actually do is that
we consider the couple consisting
of our processes issued from two different points
and define $\Theta$ and $r$ as functions of this couple.
When the starting points are close we can almost recover
the derivative of the initial
process with respect to the initial
data. Of course, the couple is a degenerate process
and that is why in \cite{Kr_1} and \cite{Kr_2}
we paid a special attention
not to impose the nondegeneracy condition whenever
it is not necessary.

In contrast with controlled diffusion processes,
no version of random time change rule, change of Wiener
process and Girsanov's theorem is known, and instead we
can only rely
on what the results of \cite{Kr_2}
allow one to extract from inspecting the corresponding
Isaacs equations.
\end{example}

%s5 #&#
%s5 ###
\section{On equivalent representations of value functions}
\label{section 5.20.1}
Here we suppose that Assumptions
\ref{assumption 5.19.1}, \ref{assumption 5.29.02},
\ref{assumption 5.19.2} and
\ref{assumption 5.22.1} are satisfied.

%
%as5.1 #&#
\begin{assumption}
\label{assumption 6.1.3}
There exists a nonnegative $
G\in C(\bar{D})\cap C^{2}_{\mathrm{loc}}(D)$
such that $G= 0$ on $\partial D$ (if $D\ne\bR^{d}$) and
\[
L^{\alpha\beta}G(p,x)\leq-1
\]
in $D$ for all $p\in\bR^{k}$, $\alpha\in A$ and $\beta\in B$.
\end{assumption}

Suppose that we are also given an $\bR^{d_{1}}$-valued
function $\pi^{\alpha\beta}(x,y)$ defined for $x,y\in\bR^{d}$,
$\alpha\in A$ and $\beta\in B$,
which is bounded by $K_{0}$, Borel measurable, and Lipschitz
continuous with respect to $x$ with Lipschitz constant $K_{1}$.

Then
for $\alpha_{\cdot}\in\frA$, $\beta_{\cdot}\in\frB$,
$x,y\in\bR^{d}$ introduce
$y^{\alpha_{\cdot}\beta_{\cdot} y}_{t}
=y^{\alpha_{\cdot}\beta_{\cdot}x,y}_{t}$
as a unique solution of the It\^o equation
%
%
%e5.1 #&#
%e5.1 ###
\begin{equation}
\label{5.7.1} y_{t}=y+\int_{0}^{t}
\sigma^{\alpha_{s}\beta_{s}} (y_{s}) \,dw_{s}+ \int
_{0}^{t}b^{\alpha_{s}\beta_{s}} (y_{s}) \,ds
\end{equation}
and introduce
$x^{\alpha_{\cdot}\beta_{\cdot}x,y}_{t}$
as a unique solution of the It\^o equation
(recall that $\hat{\sigma},\hat{b},\hat{c},\hat{f}$
are introduced before Assumption~\ref{assumption 5.19.2})
%
%
%e5.2 #&#
%e5.2 ###
\begin{equation}
\label{5.7.2} x_{t}=x+\int_{0}^{t}
\hat\sigma^{\alpha_{s}\beta_{s}} (x_{s},y_{s})
\,dw_{s} + \int_{0}^{t}( \hat b-\hat
\sigma \pi)^{\alpha_{s}\beta_{s}} (x_{s},y_{s}) \,ds,
\end{equation}
where, of course,
$y_{s}=y^{\alpha_{\cdot}\beta_{\cdot}x,y}_{s}$.
We emphasize that \eqref{5.7.1} has a unique solution
since the coefficients are Lipschitz continuous in $y$
and are bounded, and for given~$y_{\cdot}$, equation
\eqref{5.7.2} has a unique solution since its coefficients
are Lipschitz continuous in $x$ and are bounded.
It follows that, in the terminology of \cite{Kr_1}, system~\eqref{5.7.1}--\eqref{5.7.2} satisfies the usual
hypothesis [although the coefficients in
\eqref{5.7.2} may not be Lipschitz continuous
with respect to the $y$ variable].

With the above
$y_{s}$ and $x_{s}=x^{\alpha_{\cdot}\beta_{\cdot}x,y}_{s}$
also define
\[
\phi^{\alpha_{\cdot}\beta_{\cdot}x,y}_{t}= \int_{0}^{t}
\hat c^{\alpha_{s}\beta_{s}}(x_{s},y_{s}) \,ds,
\]
and for $z\in\bR$
introduce $z^{\alpha_{\cdot}\beta_{\cdot}x,y,z}_{s}$
as a unique solution of
%
%
%e5.3 #&#
%e5.3 ###
\begin{equation}
\label{5.13.1} z_{t}=z+\int_{0}^{t}z_{s}
\bigl[\pi^{\alpha_{s}\beta_{s}} (x_{s}, y_{s})
\bigr]^{*} \,dw_{s}.
\end{equation}

Next, for $X=(x,y,z)$, $x,y\in\bR^{d}$, $z\in\bR$
denote
\begin{eqnarray*}
x^{\alpha_{\cdot}\beta_{\cdot}X}_{t}&=& x^{\alpha_{\cdot}\beta_{\cdot}x,y}_{t},\qquad
y^{\alpha_{\cdot}\beta_{\cdot}X}_{t}= y^{\alpha_{\cdot}\beta_{\cdot} y}_{t},\qquad
\phi^{\alpha_{\cdot}\beta_{\cdot}X}_{t}= \phi^{\alpha_{\cdot}\beta_{\cdot}x,y}_{t}
\\
X^{\alpha_{\cdot}\beta_{\cdot}X}_{t}& =&(x_{t},y_{t},z_{t})^{\alpha_{\cdot}\beta_{\cdot}X},
\end{eqnarray*}
fix a number
$M\in(1,\infty)$, for $X=(x,y,z)$
define $\tau^{\alpha_{\cdot}\beta_{\cdot}X}$
as the first exit time of
$(x,z)^{\alpha_{\cdot}\beta_{\cdot}X}_{t}$ from $D
\times(M^{-1},M)$ and set
\[
v^{\alpha_{\cdot}\beta_{\cdot}}(X) =E^{\alpha_{\cdot}\beta_{\cdot}}_{X} \biggl[\int
_{0}^{\tau}\hat f(X_{t})e^{-\phi_{t}}
\,dt +z_{\tau}v(x_{\tau})e^{-\phi_{\tau}} \biggr],
\]
where $\hat f^{\alpha\beta}
(x,y,z)=z\hat f^{\alpha\beta}(x,y)$, and $v$ is taken
as in Theorem~\ref{theorem 5.19.1} and is at least bounded
and continuous according to the results
of \cite{Kr_2} and owing to Assumption~\ref{assumption 6.1.3}. Finally, introduce
\[
v(X)=\infsup_{\bbeta\in\mathbb{B}\,  \alpha_{\cdot}\in\frA}
v^{\alpha_{\cdot}\bbeta(\alpha_{\cdot})}(X).
\]
The fact that $v^{\alpha\beta}(X)$ and $v(X)$ are well defined
and bounded will be seen from the proof of the
following.

%
%th5.1 #&#
\begin{theorem}
\label{theorem 5.1.1}
Under the above notation for $X=(x,y,z)$
we have
%
%
%e5.4 #&#
%e5.4 ###
\begin{equation}
\label{5.15.3} v(X)=zv(x).
\end{equation}
Furthermore, if we are given stopping times $\gamma^{\alpha_{\cdot}
\beta_{\cdot}X}\leq\tau^{\alpha_{\cdot}
\beta_{\cdot}X}$, then
%
%
%e5.5 #&#
%e5.5 ###
\begin{equation}
\label{4.2.1} zv (x) =\infsup_{\bbeta\in\mathbb{B}\,  \alpha_{\cdot}\in\frA}
E_{X}^{\alpha_{\cdot}\bbeta(\alpha_{\cdot})}
\biggl[ \int_{0}^{\gamma} \hat f(X_{t})
e^{- \phi_{t} } \,dt + z_{\gamma}v (x_{\gamma}) e^{- \phi_{\gamma} }
\biggr].
\end{equation}
\end{theorem}

\begin{pf} Introduce
%
%
%e5.6 #&#
%e5.6 ###
\begin{equation}
\label{6.1.7} (\Ba,\Bsigma,\Bb,\Bc,\Bf)^{\alpha\beta}(x,y)= (a,
\sigma,b,c,f)^{\alpha\beta} \bigl( p^{\alpha\beta}(x,y),x \bigr)
\end{equation}
(specifying the value of $p$ transforms the letters to
their boldface options).
Also
denote by $\cP$ the set of triples $\check p=(r,\pi,P)$, where
$r\in[\delta_{1},\delta_{1}^{-1}]$, $\pi\in\bR^{d_{1}}$
with $|\pi|\leq K_{0}$ and $P\in\cO$.
For $\check p=(r,\pi,P)\in\cP$ define
\begin{eqnarray*}
\check\sigma^{\alpha\beta}(\check p,x,y )&=&r \Bsigma^{\alpha\beta}(x,y)P, \qquad\check
b^{\alpha\beta}(\check p,x,y )=r^{2} \Bb^{\alpha\beta}(x,y) -r
\Bsigma^{\alpha\beta}(x,y)P\pi,
\\
\check c^{\alpha\beta}(\check p,x,y,z)&=&r^{2}\Bc^{\alpha\beta}(x,y),\qquad
\check f^{\alpha\beta}(\check p,x,y,z)=r^{2}z\Bf^{\alpha\beta}(x,y)
\end{eqnarray*}
and also write
\[
r=r(\check p), \qquad\pi=\pi(\check p),\qquad P=P(\check p).
\]
We thus freed the coefficients of \eqref{5.7.2} of the particular
values of $r,\pi,P$.

For each $\check p\in\cP$ there is a natural operator $\check
L^{\alpha\beta}$
acting on smooth functions
$u(x,y,z)$ and mapping them to
\[
\check L^{\alpha\beta}u(\check p,x,y,z)
\]
associated with the matrix of second-order coefficients
\[
\frac{1}{2}\pmatrix{\check\sigma^{\alpha\beta}(\check p,x,y )\vspace*{2pt}
\cr
\sigma^{\alpha\beta}(y)\vspace*{2pt}
\cr
z \pi^{*}(\check p) }
\pmatrix{\check\sigma^{\alpha\beta}(\check p,x,y )\vspace *{2pt}
\cr
\sigma^{\alpha\beta}(y)\vspace*{2pt}
\cr
z \pi^{*}(\check p)
}^{*},
\]
the drift term
\[
\pmatrix{\check b^{\alpha\beta}(\check p,x,y )\vspace*{2pt}
\cr
\check{b}^{\alpha\beta}(y)\vspace*{2pt}
\cr
0 }
\]
and the zeroth-order
(killing) coefficient $-\check c^{\alpha\beta}(\check p,x,y,z)$.
Introduce $\bar{p}=\break (1,0,I)$ and
\[
\bar{L}^{\alpha\beta}u(x,y,z)=\check L^{\alpha\beta}u(\bar{p},x,y,z),\qquad
\bar{f}^{\alpha\beta} (x,y,z)=\check f^{\alpha\beta} (\bar{p},x,y,z).
\]
We also need the operator $\BL$ acting on functions
$u(x,y)$ by the formula
\[
\BL^{\alpha\beta}u(x,y)= \Ba^{\alpha\beta}_{ij}(x,y)D_{ij}u(x,y)+
\Bb^{\alpha\beta}_{i} (x,y)D_{i }u(x,y)-
\Bc^{\alpha\beta}(x,y) u(x,y)
\]
(no differentiation with respect to $y$ is involved).
Notice that, if $u=u(x)$ is a smooth function on $\bR^{d}$
and $\check{u}(x,y,z):=zu(x)$,
then as is easy to check
%
%
%e5.7 #&#
%e5.7 ###
\begin{equation}
\label{5.15.1} \check L^{\alpha\beta}\check{u}(\check p,x,y,z)
=zr^{2}(\check p) \bigl(\bar{L}^{\alpha\beta}u \bigr) (x,y,z)
=zr^{2}(\check p)\BL^{\alpha\beta}u(x,y).
\end{equation}

One of consequences of
Assumption~\ref{assumption 6.1.3} and \eqref{5.15.1} is that
in $D\times\bR^{d}\times(M^{-1},M)$ we have
\[
\check L^{\alpha\beta}\check{G}(\check p,x,y,z)\leq-1
\]
for all $\check p$, where $\check G(x,y,z)=M\delta_{1}^{-2}zG(x)$.
In particular, this implies that
$v^{\alpha\beta}(X)$ and $v(X)$ are well defined
and are bounded.

Next, fix $x_{0}\in D$, $y_{0}\in\bR^{d}$,
and set
\[
\check p_{t}^{\alpha_{\cdot}\beta_{\cdot}}= (r,\pi,P)^{\alpha_{t}\beta_{t}}(x,y)
^{\alpha_{\cdot}\beta_{\cdot}x_{0},y_{0}}_{t}.
\]
As is easy to see, $\check p_{t}^{\alpha_{\cdot}\beta_{\cdot}}$ is a
control adapted process in terminology of \cite{Kr_1}; see Remark~2.3 there.
For $\alpha_{\cdot}\in\frA$ and $\beta_{\cdot}\in\frB$,
consider
the following system of It\^o's equations:
%
%e5.8 ###
\begin{eqnarray}\label{5.10.3}
d\check x_{t}&=&\check \sigma^{\alpha_{t}\beta_{t}} \bigl(\check
p_{t}^{\alpha_{\cdot}\beta
_{\cdot}},\check x_{t},\check y_{t}
\bigr) \,dw_{t} +\check b^{\alpha_{t}\beta_{t}} \bigl(\check
p_{t}^{\alpha_{\cdot}\beta_{\cdot}},\check x_{t},\check y_{t}
\bigr) \,dt,\nonumber
\\
 d\check y_{t}&=&\sigma^{\alpha_{t}\beta_{t}} (\check
y_{t}) \,dw_{t} +b^{\alpha_{t}\beta_{t}}(\check y_{t})
\,dt,
\\
d\check z_{t}&=&\check z_{t} \pi^{*} \bigl(\check
p_{t}^{\alpha_{\cdot}\beta_{\cdot}} \bigr) \,dw_{t}.\nonumber
\end{eqnarray}
Its solution with initial condition
$X=(x,y,z )$
will be denoted by
\[
\check{X}^{\alpha_{\cdot}\beta_{\cdot}X}_{t} =(\check{x},\check{y},\check{z} )
^{\alpha_{\cdot}\beta_{\cdot}X}_{t}.
\]
Observe that by uniqueness,
%
%
%e5.9 #&#
%e5.9 ###
\begin{equation}
\label{7.15.1} \check{X}^{\alpha_{\cdot}\beta_{\cdot}x_{0},y_{0},z}_{t} = X^{\alpha_{\cdot}\beta_{\cdot}x_{0},y_{0},z}_{t}
\end{equation}
for any $z$.
Also define
\begin{eqnarray*}
\check{\phi}_{t}^{\alpha_{\cdot}\beta_{\cdot}X} &=&\int_{0}^{t}
\check c^{\alpha_{t}\beta_{t}} \bigl(\check p^{\alpha_{\cdot}\beta_{\cdot} }_{s},
\check{X}^{\alpha_{\cdot}\beta_{\cdot}X}_{s} \bigr) \,ds,
\\
\check{v}(X)&=&\infsup_{\bbeta\in\mathbb{B}\,  \alpha_{\cdot}\in
\frA} E_{X}^{\alpha_{\cdot}\bbeta(\alpha_{\cdot})} \biggl[
\int_{0}^{\check{\tau}} \check f(\check p_{t},
\check X_{t})e^{-\check{\phi}_{t}} \,dt+ \check{z}_{\check{\tau}}v(
\check{x}_{\check{\tau}}) e^{-\check{\phi}_{\check{\tau}}} \biggr],
\end{eqnarray*}
where $\check{\tau}^{\alpha_{\cdot}\beta_{\cdot}X}$
is the first exit time of
$\check{X}^{\alpha_{\cdot}\beta_{\cdot}X}_{t}$
from $D{\V}=D\times\bR^{d}\times(M^{-1},M)$.

It turns out that, in the terminology of \cite{Kr_1},
for any $C^{2}_{\mathrm{loc}}(D)$ function $u=u(x)$,
the function $zu(x)$ is $p$-insensitive
in $D{\V}$ relative
to
$( zr^{2}(\check p),
\check L^{\alpha\beta})$. This follows from the fact that,
if $ X
\in D{\V}$, then
by It\^o's formula and \eqref{5.15.1}, for
$t<\check{\tau}^{\alpha_{\cdot}\beta_{\cdot}X}$,
\begin{eqnarray*}
&&d \bigl(u \bigl(\check{x}^{\alpha_{\cdot}\beta_{\cdot}X} _{t} \bigr)
\check{z}^{\alpha_{\cdot}\beta_{\cdot}X}_{t} e^{-\check{\phi}^{\alpha_{\cdot}\beta_{\cdot}X}_{t}} \bigr)
\\
&&\qquad=e^{-\check{\phi}_{t}} \check{z}^{\alpha_{\cdot}\beta_{\cdot}X}_{t} r^{2} \bigl(
\check p^{\alpha_{\cdot}\beta_{\cdot}}_{t} \bigr) \bigl(\bar{L}^{\alpha_{t}\beta_{t}}u
\bigr) \bigl(\check{x}^{\alpha_{\cdot}\beta_{\cdot}X}_{t}, \check y_{t}^{\alpha_{\cdot}\beta_{\cdot}X},
\check z_{t}^{\alpha_{\cdot}\beta_{\cdot}X} \bigr) \,dt+dm_{t},
\end{eqnarray*}
where $m_{t}$ is a local martingale starting at zero,
and $zr^{2}(\check p)\in[M^{-1}\delta_{1}^{2},
M\delta_{1}^{-2}]$.

Furthermore, it turns out that
equation \eqref{5.15.1} and Assumption~\ref{assumption 5.22.1}
also imply that for smooth $u=u(x)$, if
at a particular point $x$ it holds that
\[
J(x):=\supinf_{\alpha\in A\,  \beta\in B} \bigl[ a^{\alpha\beta}_{ij}(x)D_{ij}u(x
)+ b^{\alpha\beta}_{i} (x)D_{i }u(x )- c^{\alpha\beta}(x)
u(x )+ f^{\alpha\beta} (x ) \bigr]\leq0,
\]
then with the same $x$, any $y$ and $z>0$, we also have
\[
I(x,y,z):=\supinf_{\alpha\in A\,  \beta\in B} \bigl[\bar{L}^{\alpha\beta}\check{u}(x,y,z)+
\bar{f}^{\alpha\beta} (x,y,z) \bigr] \leq0,
\]
where $\check{u}(x,y,z):=zu(x)$.
Indeed, since
\begin{eqnarray*}
&&J(x)=\supinf_{\alpha\in A\,  \beta\in B} \bigl[ \Ba^{\alpha\beta}_{ij}(x,x)D_{ij}u(x
)+ \Bb^{\alpha\beta}_{i} (x,x)D_{i }u(x )
- \Bc^{\alpha\beta}(x,x) u(x )\\
&&\hspace*{260pt}{}+ \Bf^{\alpha\beta} (x,x ) \bigr],
\end{eqnarray*}
the inequality $J(x)\leq0$ implies
by Assumption~\ref{assumption 5.22.1}
that
\begin{eqnarray*}
&&\supinf_{\alpha\in A\,  \beta\in B} \bigl[ \Ba^{\alpha\beta}_{ij}(x,y)D_{ij}u(x
)+ \Bb^{\alpha\beta}_{i} (x,y)D_{i }u(x )
- \Bc^{\alpha\beta}(x,y) u(x )\\
&&\hspace*{224pt}{}+ \Bf^{\alpha\beta} (x,y ) \bigr]\leq0,
\end{eqnarray*}
and it only remains to notice that
the left-hand side is just $z^{-1}I(x,y,z)$.
Similarly, $J(x)\geq0$ implies that $I(x,y,z)\geq0$.

These facts combined imply by Theorems 2.3
and 3.1 of \cite{Kr_2}
that
for all $x\in\bar{D}$, $y\in\bR^{d}$ and $z\in
[M^{-1},M]$ we have
\[
\check{v}(x,y,z)=zv(x)
\]
and, for any stopping
times $\gamma^{\alpha_{\cdot}\beta_{\cdot}X}
\leq\check{\tau}^{\alpha_{\cdot}\beta_{\cdot}X}$,
%
%
%e5.10 #&#
%e5.10 ###
\begin{equation}
\label{7.15.2} zv (x) =\infsup_{\bbeta\in\mathbb{B}\,  \alpha_{\cdot}\in\frA} E_{X}^{\alpha_{\cdot}\bbeta(\alpha_{\cdot})}
\biggl[ \int_{0}^{\gamma} \check f(\check
p_{t},\check X_{t}) e^{-\check\phi_{t} } \,dt + \check
z_{\gamma}v (\check x_{\gamma}) e^{- \check\phi_{\gamma} } \biggr].
\end{equation}

By \eqref{7.15.1} for $X_{0}=(x_{0},y_{0},
z_{0})$, $z_{0}\in
[M^{-1},M]$, we have
\begin{eqnarray*}
\check X^{\alpha_{\cdot}\beta_{\cdot}X_{0}} &=&X^{\alpha_{\cdot}
\beta_{\cdot}X_{0}},\qquad \check f^{\alpha_{t}\beta_{t}} \bigl(
\check{p}^{\alpha_{\cdot}
\beta_{\cdot}}_{t},\check X_{t} ^{\alpha_{\cdot}\beta_{\cdot}X_{0}}
\bigr)=\hat f \bigl(X_{t} ^{\alpha_{\cdot}\beta_{\cdot}X_{0}} \bigr),
\\
\check\phi^{\alpha_{\cdot}\beta_{\cdot}X_{0}}& =&\phi^{\alpha_{\cdot}\beta_{\cdot}X_{0}},
\end{eqnarray*}
so that
$v(x_{0},y_{0},
z_{0})=\check{v}(x_{0},y_{0},
z_{0})$. It follows that \eqref{5.15.3} holds at
$(x_{0},y_{0},
z_{0})\in D{\V}$. Outside
$D{\V}$ the equality is obvious. Finally,
\eqref{4.2.1} follows from \eqref{7.15.2},
and the theorem is proved.
\end{pf}

%
%re5.1 #&#
\begin{remark}
One of assumptions in Theorems 2.3 and 3.1 of \cite{Kr_2}
is that the coefficients satisfy
Assumption~\ref{assumption 5.19.1}(i) without $p^{\alpha\beta}(x,y)$ there. Since $p$ is involved in
\eqref{6.1.7} we needed to include it in
Assumption~\ref{assumption 5.19.1}(i) in contrast with the parameters
$r^{\alpha\beta}(x,y)$ and $P^{\alpha\beta}(x,y)$.
The same reasons caused the last requirement in
Assumption~\ref{assumption 5.19.1}(ii). Recall that in
Theorems 2.3 and 3.1 of \cite{Kr_2} the coefficients of
It\^o equations are not supposed to be Lipschitz,
but rather uniformly continuous.
\end{remark}

%s6 #&#
%s6 ###
\section{Estimating the difference of solutions of stochastic
equations whose coefficients are close}

\label{section 5.20.2}

Suppose that on $\Omega\times(0,\infty)\times\bR^{d}
\times\bR^{d}$
we are given the following functions: $d\times d_{1}$
matrix-valued $\sigma_{t}(x,y) $,
$\bR^{d}$-valued $b_{t}(x,y)$ and
real-valued functions $c_{t}(x,y)\geq\delta_{1}$,
$f_{t}(x,y) $, where $\delta_{1}>0$ is a fixed
constant.

Introduce
\[
(\sigma_{t},b_{t},c_{t},f_{t}) (x)=
(\sigma_{t},b_{t}, c_{t}, f_{t})
(x,x), \qquad a_{t}(x)=(1/2)\sigma_{t}\sigma^{*}_{t}(x).
\]

%
%as6.1 #&#
\begin{assumption}
\label{assumption 5.17.1}
(i) All the above functions are measurable with respect to the
product of $\cF$ and Borel $\sigma$-algebras on $(0,\infty)$,
$\bR^{d}$ and $\bR^{d}$, and they are
progressively measurable as functions of $(\omega,t)$
for each $(x,y)$.\vspace*{-6pt}
\begin{longlist}[(iii)]
\item[(ii)] All the above functions are bounded by a constant $K_{0}$.

\item[(iii)] For any $t>0$, $x',x'',y\in\bR^{d}$
and
\[
\xi_{t}=(\sigma_{t},b_{t} ) (x,y),\qquad
\eta_{t}=(\sigma_{t},b_{t} ) (x),
\]
we have
\[
\bigl|\xi_{t} \bigl(x',y \bigr)-\xi_{t}
\bigl(x'',y \bigr)\bigr| +\bigl|\eta_{t}
\bigl(x' \bigr)-\eta_{t} \bigl(x ''
\bigr)\bigr|\leq K_{1}\bigl|x'-x''\bigr|,
\]
where $K_{1}$ is a fixed constant. Also there exists
a constant $\varepsilon_{0}>0$ such that
for any $t>0$ and $x,y\in\bR^{d}$ with $|x-y|
\leq\varepsilon_{0}$, we have
\[
\bigl|c_{t}(x,y)-c_{t}(y)\bigr|+\bigl|f_{t}(x,y)-f_{t}(y)\bigr|
\leq K_{1}|x-y|.
\]
\end{longlist}
\end{assumption}

Observe that Assumption~\ref{assumption 5.17.1}(iii)
implies, in particular, that $|b_{t}(x,y)-b_{t}(y)|\leq
K_{1}|x-y|$.\vadjust{\goodbreak}

%
%as6.2 #&#
\begin{assumption}
\label{assumption 4.24.1}
There exist constants $\mu\geq1$
and $\delta\geq2\delta_{1} $ such that
for all $x,y\in\bR^{d}$ satisfying
$|x-y|\leq\varepsilon_{0}$ we have
%
%e6.1 ###
\begin{eqnarray}\label{4.24.2}
R_{t}(x,y)&:=&\bigl\|\sigma_{t} (x, y)-\sigma_{t} (y)
\bigr\|^{2}_{\xi} + 2 \bigl\langle x-y,b_{t} (x, y)
-b_{t}(y) \bigr\rangle
\nonumber
\\
 &&{}-4\mu \bigl\langle x-y, a_{t}(x) (x-y) \bigr\rangle
\\
&\leq&2 \bigl(c_{t}(y) -\delta \bigr)|x-y|^{2},\nonumber
\end{eqnarray}
where $\xi=(x-y)/|x-y|$.
\end{assumption}

Fix a unit $\xi\in\bR^{d}$, and for $\varepsilon\in
[0,\varepsilon_{0}]$
introduce $x_{t}^{\varepsilon}$ as a unique solution of
\[
x_{t}=\varepsilon\xi+\int_{0}^{t}
\sigma_{s} (x_{s},y_{s}) \,dw_{s}
+ \int_{0}^{t} \bigl[b_{s}
(x_{s},y_{s})- 2\mu a_{s}(x_{s})
(x_{s}-y_{s}) \bigr] \,ds,
\]
where $y_{s} $ is a unique solution of
\[
y_{t}= \int_{0}^{t}
\sigma_{s}(y_{s}) \,dw_{s} +\int
_{0}^{t}b_{s}(y_{s}) \,ds.
\]
Observe that owing to uniqueness,
\[
x_{t}^{0}=y_{t}.
\]

For $\varepsilon>0$ define
\[
\xi^{\varepsilon}_{t}=\frac{1}{\varepsilon} \bigl(x^{\varepsilon}_{t}-x^{0}_{t}
\bigr),\qquad \phi_{t}=\int_{0}^{t}c_{s}
\bigl(x^{0}_{s} \bigr) \,ds,
\]
and for $\lambda>0$ let
\[
\kappa_{\varepsilon}(\lambda) =\inf \bigl\{t\geq0\dvtx
 \bigl|x^{\varepsilon}_{t}-x^{0}_{t}\bigr|
\geq\lambda \bigr\}.
\]

Notice that $\kappa_{\varepsilon}(\lambda)=0$
if $\lambda\leq\varepsilon$, and start with the following:

%
%le6.1 #&#
\begin{lemma}
\label{lemma 2.3.1}
For any $\lambda\in(0,\varepsilon_{0}]$
%
%
%e6.2 #&#
%e6.3 ###
%e6.2 ###
\begin{eqnarray}
\label{4.7.1} J_{\varepsilon}&:=&E\int_{0}^{\kappa_{\varepsilon}(\lambda)} \bigl|
\xi_{t}^{\varepsilon}\bigr|e^{-\phi_{t}
+\delta t /2} \,dt\leq2/\delta,
\\
\label{2.6.1} I_{\varepsilon}&:=&E\sup_{t<\kappa_{\varepsilon}(\lambda)}\bigl|
\xi_{t}^{\varepsilon}\bigr|e^{-\phi_{t}+\delta t /2}\leq N,
\end{eqnarray}
where $N$ is a constant depending only on $K_{1}$
and $\delta$.
\end{lemma}

\begin{pf} We have
%
%e6.4 ###
\begin{eqnarray}\label{4.7.3}
d\xi_{t}^{\varepsilon} &=& \varepsilon^{-1} \bigl[
\sigma_{t} \bigl(x^{\varepsilon}_{t},x^{0}_{t}
\bigr)-\sigma_{t} \bigl(x^{0}_{t} \bigr) \bigr]
\,dw_{t}
\nonumber
\\[-8pt]
\\[-8pt]
\nonumber
&&{} +\varepsilon^{-1} \bigl[ b_{t}
\bigl(x^{\varepsilon}_{t},x^{0}_{t} \bigr) -b
\bigl(x^{0}_{t} \bigr)-2\mu a_{t}
\bigl(x^{\varepsilon}_{t} \bigr) \bigl(x^{\varepsilon}_{t}-x^{0}_{t}
\bigr) \bigr] \,dt,
\end{eqnarray}
where the magnitudes of the coefficients of $dw_{t}$ and $dt$ are dominated
by constants times $|\xi^{\varepsilon}_{t}|$. This
allows us to use It\^o's formula (cf. the proof of Theorem~5.8.7 of~\cite{Intro}) and obtain that ($0/0:=0$)
\begin{eqnarray*}
&&d\bigl|\varepsilon\xi_{t}^{\varepsilon}\bigr|e^{-\phi_{t}+\delta t /2}\\
&&\qquad=
\frac{1}{2|x^{\varepsilon}_{t}-x^{0}_{t}|} \bigl[R_{t}
\bigl( x^{\varepsilon}_{t},
x^{0}_{t} \bigr)-2 \bigl(c_{t}
\bigl(x^{0}_{t} \bigr)-\delta/2 \bigr)\bigl|x^{\varepsilon}_{t}-
x^{0}_{t}\bigr|^{2} \bigr] e^{-\phi_{t}+\delta t /2} \,dt
\\
&&\qquad\quad{}+S_{t} \bigl( x^{\varepsilon}_{t}, x^{0}_{t}
\bigr)e^{-\phi_{t}+\delta t /2} \,dw_{t},
\end{eqnarray*}
where
\[
S_{t} \bigl( x^{\varepsilon}_{t}, x^{0}_{t}
\bigr)=\frac{1}{| \xi_{t}^{\varepsilon}|} \xi_{t}^{\varepsilon*} \bigl[
\sigma_{t} \bigl(x^{\varepsilon}_{t},x^{0}_{t}
\bigr)- \sigma_{t} \bigl(x^{0}_{t} \bigr)
\bigr].
\]
By assumption, for
$t<\kappa_{\varepsilon}(\lambda)$ we have
\[
R_{t} \bigl( x^{\varepsilon}_{t}, x^{0}_{t}
\bigr)-2 \bigl(c_{t} \bigl(x^{0}_{t} \bigr)-
\delta/2 \bigr)\bigl|x^{\varepsilon}_{t}- x^{0}_{t}\bigr|^{2}
\leq- \delta\bigl|x^{\varepsilon}_{t}- x^{0}_{t}\bigr|^{2}.
\]
It follows that for
$t<\kappa_{\varepsilon}(\lambda)$,
%
%
%e6.5 #&#
%e6.5 ###
\begin{equation}
\label{6.4.1}\qquad d\bigl| \xi_{t}^{\varepsilon}\bigr|e^{-\phi_{t}+\delta t /2} \leq- (
\delta/2) \bigl| \xi_{t}^{\varepsilon}\bigr|e^{-\phi_{t}+\delta t
/2} \,dt +
\varepsilon^{-1}S_{t} \bigl( x^{\varepsilon}_{t},
x^{0}_{t} \bigr)e^{-\phi_{t}+\delta t /2} \,dw_{t}.
\end{equation}
In particular, \eqref{4.7.1} holds.
Furthermore,
%
%
%e6.6 #&#
%e6.6 ###
\begin{equation}
\label{6.4.2} \bigl|\varepsilon^{-1}S_{t} \bigl(
x^{\varepsilon}_{t}, x^{0}_{t} \bigr)\bigr|\leq
K_{1}\bigl|\xi^{\varepsilon}_{t}\bigr|,
\end{equation}
and by Davis's inequality,
\begin{eqnarray*}
I_{\varepsilon} &\leq&3K_{1}E \biggl(\int_{0}^{\kappa_{\varepsilon}(\lambda)}
\bigl| \xi_{t}^{\varepsilon}\bigr|^{2}e^{-2\phi_{t}+\delta t } \,dt
\biggr)^{1/2}
\\
&\leq&3K_{1}E \Bigl(\sup_{s<\kappa_{\varepsilon}(\lambda)}\bigl|
\xi_{s}^{\varepsilon}\bigr|e^{-\phi_{s}+\delta s /2} \Bigr)^{1/2} \biggl(
\int_{0}^{\kappa_{\varepsilon}(\lambda)} \bigl| \xi_{t}^{\varepsilon}\bigr|e^{-\phi_{t}+\delta t /2}
\,dt \biggr)^{1/2}\leq NI^{1/2}_{\varepsilon}J^{1/2}_{\varepsilon},
\end{eqnarray*}
which, due to \eqref{4.7.1}, proves \eqref{2.6.1}
and the lemma.
\end{pf}

%
%co6.2 #&#
\begin{corollary}
\label{corollary 4.7.1}
For $\lambda>0$ we have
\[
Ee^{-\phi_{\kappa_{\varepsilon}(\lambda)}+\kappa_{\varepsilon}
(\lambda)\delta/2}I_{\kappa_{\varepsilon}(\lambda)<\infty} \leq N\varepsilon/\lambda.
\]
\end{corollary}

Indeed, if $\lambda\leq\varepsilon$, the estimate
is obvious since $\kappa_{\varepsilon}(\lambda)=0$
and for $\lambda>\varepsilon$
\[
\lambda Ee^{-\phi_{\kappa_{\varepsilon}(\lambda)}
+\kappa_{\varepsilon}(\lambda)\delta/2} I_{\kappa_{\varepsilon}(\lambda)<\infty}
=\varepsilon E\bigl|\xi^{\varepsilon}_{\kappa_{\varepsilon}(\lambda)}\bigr|
e^{-\phi_{\kappa_{\varepsilon}(\lambda)}
+\kappa_{\varepsilon}(\lambda)\delta/2}
I_{\kappa_{\varepsilon}(\lambda)<\infty} \leq N\varepsilon.
\]
%
%
%re6.1 #&#
\begin{remark}
\label{remark 6.4.1}
If $\delta\geq K_{1}^{2}$, then it follows from
\eqref{6.4.1} and \eqref{6.4.2} that for $t<\kappa_{\varepsilon}
(\lambda)$ we have
\[
d\bigl| \xi_{t}^{\varepsilon}\bigr|^{2}e^{-2\phi_{t}+\delta t } \leq
dm_{t},
\]
where $m_{t}$ is a local martingale. Hence,
for any stopping time $\gamma\leq\kappa_{\varepsilon}
(\lambda)$,
\[
E\bigl| \xi_{\gamma}^{\varepsilon}\bigr|^{2} e^{-2\phi_{\gamma}+\delta\gamma}\leq1.
\]
\end{remark}
Psychologically, the condition $\delta\geq K_{1}^{2}$
may look artificial. However, in the proof of Theorem~\ref{theorem 6.4.1} the parameter $\delta$ will be, basically,
sent to infinity.

Next introduce
\[
\pi_{s}(x,y)=\mu\sigma_{s}^{*}(x) (x-y)
\]
and
introduce $\rho_{t}^{\varepsilon} $ as a unique
solution of
\[
\rho_{t}=1+\int_{0}^{t}
\rho_{s}\pi^{*} _{s} \bigl(x^{\varepsilon}_{s},x^{0}_{s}
\bigr) \,dw_{s}+\int_{0}^{t}
\rho_{s} \bigl[c_{s} \bigl(x^{0}_{s}
\bigr)- c_{s} \bigl(x^{\varepsilon}_{s},x^{0}_{s}
\bigr) \bigr] \,ds.
\]
Take a constant $M>1$ and define
\[
\gamma_{\varepsilon}(M)
\]
as the first exit time of $\rho^{\varepsilon}_{t} $
from $(M^{-1},M)$.

Recall that $c\geq\delta_{1}$.
%
%
%le6.3 #&#
\begin{lemma}
\label{lemma 4.7.3}
There exists $\lambda_{1}\in(0,\varepsilon_{0}]$,
depending only on $\varepsilon_{0},K_{0},K_{1}$
and $\delta_{1}$,
and there exists a constant $N$,
depending only on $ K_{1}$
and $\delta_{1}$, such that for $\lambda=\lambda_{1}/\mu$
and $\mu\geq1$ we have
%
%
%e6.7 #&#
%e6.7 ###
\begin{equation}
\label{4.7.5} I:=E\sup_{t<\gamma_{\varepsilon}(M)\wedge\kappa_{\varepsilon
}(\lambda)} \bigl|\rho_{t}^{\varepsilon}-1\bigr|
e^{-\phi_{t}+\delta_{1}t/2 } \leq N \bigl(M\mu^{2}+1 \bigr)^{1/2}
\delta^{-1/2}\varepsilon.
\end{equation}
\end{lemma}

\begin{pf} Denote $C_{t}(x^{\varepsilon}_{t},
x^{0}_{t})=c_{t}(x^{0}_{t})-c_{t}(x^{\varepsilon}_{t},
x^{0}_{t})$ and
$
\eta_{t}=(\rho^{\varepsilon}_{t}-1)^{2}
$.
Then
\begin{eqnarray*}
d\eta_{t}&=&2 \bigl(\rho^{\varepsilon}_{t}-1 \bigr)
\rho^{\varepsilon}_{t}\pi^{*}_{t}
\bigl(x^{\varepsilon}_{t}, x^{0}_{t} \bigr)
\,dw_{t} +2 \bigl(\rho^{\varepsilon}_{t}-1 \bigr)
\rho^{\varepsilon}_{t} C_{t} \bigl(x^{\varepsilon}_{t},
x^{0}_{t} \bigr) \,dt \\
&&{}+\bigl|\rho^{\varepsilon}_{t}\bigr|^{2}
\bigl|\pi_{t} \bigl(x^{\varepsilon}_{t},x^{0}_{t}
\bigr)\bigr|^{2} \,dt,
\\
d\eta_{t}e^{-2\phi_{t} +\delta_{1}t }&=& e^{-2\phi_{t} +\delta_{1}t } \bigl[ 2\eta_{t}
C_{t} \bigl(x^{\varepsilon}_{t}, x^{0}_{t}
\bigr)+2 \bigl(\rho^{\varepsilon}_{t}-1 \bigr) C_{t}
\bigl(x^{\varepsilon}_{t}, x^{0}_{t} \bigr)
\\
&&\hspace*{44pt}{}+\eta_{t}\bigl|\pi_{t} \bigl(x^{\varepsilon}_{t},x^{0}_{t}
\bigr)\bigr|^{2} +2 \bigl(\rho^{\varepsilon}_{t}-1 \bigr)\bigl|
\pi_{t} \bigl(x^{\varepsilon}_{t},x^{0}_{t}
\bigr)\bigr|^{2}\\
&&\hspace*{58pt}\qquad{} +\bigl|\pi_{t} \bigl(x^{\varepsilon}_{t},x^{0}_{t}
\bigr)\bigr|^{2} -\eta_{t} \bigl(2c_{t}
\bigl(x^{0}_{t} \bigr)-\delta_{1} \bigr) \bigr]
\,dt+dm_{t},
\end{eqnarray*}
where $m_{t}$ is a local martingale starting at zero,
and for $t<\gamma_{\varepsilon}(M)$,
the expression in the square brackets is less than
\begin{eqnarray*}
&&\eta_{t} \bigl[2C_{t} \bigl(x^{\varepsilon}_{t},
x^{0}_{t} \bigr)+ \delta_{1}/2+\bigl|
\pi_{t} \bigl(x^{\varepsilon}_{t},x^{0}_{t}
\bigr)\bigr|^{2} - \bigl(2c_{t} \bigl(x^{0}_{t}
\bigr) -\delta_{1} \bigr) \bigr]
\\
&&\qquad{}+(2/\delta_{1})C^{2}_{t} \bigl(x^{\varepsilon}_{t},
x^{0}_{t} \bigr) +(2M-1)\bigl |\pi_{t}
\bigl(x^{\varepsilon}_{t},x^{0}_{t}
\bigr)\bigr|^{2}.
\end{eqnarray*}
We have that $ |G_{t}|\leq K_{1}
| x^{\varepsilon}_{t}-
x^{0}_{t} |$, $ |\pi_{t}|\leq\mu K_{0}| x^{\varepsilon}_{t}-
x^{0}_{t} |$, $c\geq\delta_{1}$ and $\mu\geq1$ and, therefore,
one can find $\lambda_{1}\in(0,\varepsilon_{0}]$
such that, for $\lambda=\lambda_{1}/\mu$
and $t<\kappa_{\varepsilon}(\lambda)$,
\[
2C_{t} \bigl(x^{\varepsilon}_{t}, x^{0}_{t}
\bigr)+ \delta_{1}/2+\bigl|\pi_{t} \bigl(x^{\varepsilon}_{t},x^{0}_{t}
\bigr)\bigr|^{2} - \bigl(2c_{t} \bigl(x^{0}_{t}
\bigr) -\delta_{1} \bigr)\leq0
\]
and then
\[
d\eta_{t} e^{-2\phi_{t} +\delta_{1}t } \leq N_{1} \bigl(M
\mu^{2}+1 \bigr) \varepsilon^{2}\bigl|\xi^{\varepsilon}_{t}\bigr|^{2}
e^{-2\phi_{t} +\delta_{1}t} \,dt+dm_{t}.
\]

Hence, for any bounded stopping time $\tau$ it holds that
\begin{eqnarray*}
&&E\eta_{\tau\wedge\gamma_{\varepsilon}(M)\wedge
\kappa_{\varepsilon}(\lambda)} e^{-2\phi_{\tau\wedge\gamma_{\varepsilon}(M)\wedge
\kappa_{\varepsilon}(\lambda)}
+\delta_{1}(\tau\wedge\gamma_{\varepsilon}(M)\wedge
\kappa_{\varepsilon}(\lambda) ) }
\\
&&\qquad\leq N_{1} \bigl(M\mu^{2}+1 \bigr) \varepsilon^{2}
E\int_{0}^{\tau\wedge\gamma_{\varepsilon}(M)\wedge
\kappa_{\varepsilon}(\lambda)} \bigl|\xi^{\varepsilon}_{t}\bigr|^{2}
e^{-2\phi_{t} +\delta_{1}t} \,dt,
\end{eqnarray*}
which owing to well-known properties of such inequalities
(see, e.g., Theorem~3.6.8 in \cite{Intro}) implies that
\begin{eqnarray*}
&& E\sup_{t\leq\gamma_{\varepsilon}(M)
\wedge\kappa_{\varepsilon}(\lambda)}\eta_{t}^{1/2}
e^{- \phi_{t} +\delta_{1}t/2 }
\\
&&\qquad\leq3N_{1} \bigl(M\mu^{2}+1 \bigr)^{1/2}
\varepsilon E \biggl(\int_{0}^{ \kappa_{\varepsilon}(\lambda)} \bigl|
\xi^{\varepsilon}_{t}\bigr|^{2} e^{-2\phi_{t} +\delta_{1}t } \,dt
\biggr)^{1/2}.
\end{eqnarray*}
Owing to \eqref{2.6.1} and the assumption that
$\delta\geq2\delta_{1}$, the last expectation
is dominated by
\[
N \biggl(\int_{0}^{\infty}e^{(\delta_{1}-\delta)t} \,dt
\biggr)^{1/2} \leq N\delta^{-1/2}.
\]
The lemma is proved.
\end{pf}

%
%co6.4 #&#
\begin{corollary}
\label{corollary 2.14.1}
There is a constant $N$,
depending only on $ K_{1}$
and $\delta_{1}$, such that for any $M\geq2$
and $ \lambda=\lambda_{1}/\mu$
%
%
%e6.8 #&#
%e6.8 ###
\begin{equation}
\label{2.14.2} Ee^{-\phi_{\gamma_{\varepsilon}(M)
\wedge\kappa_{\varepsilon}(\lambda)}} \leq N \bigl[ \mu+ \bigl(M\mu^{2}+1
\bigr)^{1/2}\delta^{-1/2} \bigr]\varepsilon.
\end{equation}
\end{corollary}

To prove \eqref{2.14.2}, it suffices
to notice that
\begin{eqnarray*}
Ee^{-\phi_{\gamma_{\varepsilon}(M)
\wedge\kappa_{\varepsilon}(\lambda)}} I_{\gamma_{\varepsilon}(M)
<\kappa_{\varepsilon}(\lambda)} &\leq& M(M-1)^{-1}E\bigl|
\rho^{\varepsilon}_{\gamma_{\varepsilon}(M)}-1\bigr| e^{-\phi_{\gamma_{\varepsilon}(M)
}}I_{\gamma_{\varepsilon}(M)
<\kappa_{\varepsilon}}
\\
&\leq& M(M-1)^{-1}E \sup_{t<\gamma_{\varepsilon}(M)
\wedge\kappa_{\varepsilon}(\lambda)} \bigl|\rho^{\varepsilon}_ {t }-1\bigr|
e^{-\phi_{t} }
\end{eqnarray*}
and then to use Corollary~\ref{corollary 4.7.1} and
to recall that $c\geq\delta_{1}$.

Now for $\lambda=\lambda_{1}/\mu$,
$\varepsilon\in(0,\varepsilon_{0}]$,
and $M\geq2$ take a stopping time
\[
\tau\leq\gamma_{\varepsilon}(M)\wedge\kappa_{\varepsilon
}(\lambda).
\]
Also take a function $g_{t}(x)$, which is measurable
in $(\omega,t,x)$ and such that $|g|\leq K_{0}$ and introduce
\[
v^{\varepsilon}= E \biggl[\int_{0}^{\tau}
z^{\varepsilon}_{t} f \bigl(x^{\varepsilon}_{t},x^{0}_{t}
\bigr)e^{-\phi^{\varepsilon}_{t}} \,dt + z^{\varepsilon}_{\tau} g_{\tau}
\bigl(x^{\varepsilon}_{\tau} \bigr) e^{-\phi^{\varepsilon}_{\tau}} \biggr],\vadjust{\goodbreak}
\]
where
\[
\phi^{\varepsilon}_{t}=\int_{0}^{t}c_{s}
\bigl( x^{\varepsilon}_{s},x^{0}_{s} \bigr)
\,ds
\]
and $ z^{\varepsilon}_{t}$ is defined as a unique solution of
\[
z_{t}=1+\int_{0}^{t}
z_{s}\pi^{*} _{s} \bigl(x^{\varepsilon}_{s},x^{0}_{s}
\bigr) \,dw_{s}.
\]
Finally, define
\[
v^{0}= E \biggl[\int_{0}^{\tau} f
\bigl(x^{0}_{t} \bigr)e^{-\phi_{t}} \,dt +
g_{\tau} \bigl(x^{0}_{\tau} \bigr) e^{-\phi_{\tau}}
\biggr].
\]

%
%th6.5 #&#
\begin{theorem}
\label{theorem 4.6.1}
Suppose that there is a constant $N_{0}$ such that
%
%
%e6.9 #&#
%e6.9 ###
\begin{equation}
\label{5.18.1} E\bigl|g_{\tau} \bigl(x^{\varepsilon}_{\tau}
\bigr)-g_{\tau} \bigl(x^{0}_{\tau}
\bigr)\bigr|e^{-\phi
_{\tau}} I_{\tau<\gamma_{\varepsilon}(M)
\wedge\kappa_{\varepsilon}(\lambda)} \leq N_{0}\varepsilon.
\end{equation}
Then there exists a constant $N$, depending only on
$K_{0}$, $K_{1}$ and $\delta_{1}$, such that for $
\lambda=\lambda_{1}/\mu$ we have
\[
\bigl|v^{\varepsilon}-v^{0}\bigr|\leq N_{0}\varepsilon+N \bigl[\mu+
\bigl(M\mu^{2}+1 \bigr)^{1/2} \delta^{-1/2}+
\delta^{-1} \bigr]\varepsilon.
\]
\end{theorem}

\begin{pf} First notice that
\[
z^{\varepsilon}_{t}e^{-\phi^{\varepsilon}_{t}}= \rho^{\varepsilon}_{t}e^{-\phi_{t}},
\]
so that
\[
\biggl|\int_{0}^{\tau} \bigl[ z^{\varepsilon}_{t}
f \bigl(x^{\varepsilon}_{t},x^{0}_{t}
\bigr)e^{-\phi^{\varepsilon}_{t}} - f \bigl(x^{0}_{t}
\bigr)e^{-\phi_{t}} \bigr] \,dt \biggr|\leq I_{\varepsilon}+J_{\varepsilon},
\]
where
\begin{eqnarray*}
I_{\varepsilon}&=&\int_{0}^{\tau}\bigl|
\rho^{\varepsilon}_{t}-1\bigr| \bigl|f \bigl(x^{\varepsilon}_{t},x^{0}_{t}
\bigr)\bigr|e^{-\phi_{t}} \,dt,
\\
J_{\varepsilon}&=&\int_{0}^{\tau}\bigl| f
\bigl(x^{\varepsilon}_{t},x^{0}_{t} \bigr) - f
\bigl(x^{0}_{t} \bigr) \bigr|e^{-\phi_{t}} \,dt.
\end{eqnarray*}

By Lemma~\ref{lemma 4.7.3},
\begin{eqnarray*}
EI_{\varepsilon}&\leq& NE\sup_{s\leq\tau}\bigl|\rho^{\varepsilon}_{s}-1\bigr|
e^{-\phi_{s}+\delta_{1} s/2} \int_{0}^{\infty}e^{-\delta_{1} t/2}
\,dt
\\
&\leq& N \bigl(M\mu^{2}+1 \bigr)^{1/2} \delta^{-1/2}
\varepsilon.
\end{eqnarray*}

By Lemma~\ref{lemma 2.3.1},
\[
EJ_{\varepsilon}\leq N\varepsilon E\int_{0}^{\tau}
\bigl|\xi^{\varepsilon}_{t}\bigr|e^{-\phi_{t}} \,dt\leq N\varepsilon/\delta.\vadjust{\goodbreak}
\]

Next
\begin{eqnarray*}
E\bigl| z^{\varepsilon}_{\tau} g_{\tau} \bigl(x^{\varepsilon}_{\tau}
\bigr) e^{-\phi^{\varepsilon}
_{\tau}} -g_{\tau} \bigl(x^{0}_{\tau}
\bigr) e^{-\phi
_{\tau}}\bigr|&=&E\bigl|\rho^{\varepsilon}_{\tau} g_{\tau}
\bigl(x^{\varepsilon}_{\tau} \bigr) -g_{\tau}
\bigl(x^{0}_{\tau} \bigr)\bigr| e^{-\phi
_{\tau}}
\\
&\leq& K_{0}E\bigl|\rho_{\tau}^{\varepsilon}-1\bigr|e^{-\phi_{\tau}} +E\bigl|
g_{\tau} \bigl(x^{\varepsilon}_{\tau} \bigr) -g_{\tau}
\bigl(x^{0}_{\tau} \bigr)\bigr| e^{-\phi
_{\tau}},
\end{eqnarray*}
where the first term is estimated as above and, owing to
\eqref{5.18.1}, the second term is dominated by
\begin{eqnarray*}
&&N_{0}\varepsilon+E\bigl| g_{\tau} \bigl(x^{\varepsilon}_{\tau}
\bigr) -g_{\tau} \bigl(x^{0}_{\tau} \bigr)\bigr|
e^{-\phi
_{\tau}} I_{\tau=\gamma_{\varepsilon}(M)
\wedge\kappa_{\varepsilon}(\lambda)}
\\
&&\qquad\leq N_{0}\varepsilon+ 2K_{0}Ee^{-\phi
_{\gamma_{\varepsilon}(M)
\wedge\kappa_{\varepsilon}(\lambda)}}\leq
N_{0}\varepsilon+N \bigl[ \mu+ \bigl(M\mu^{2}+1
\bigr)^{1/2}\delta^{-1/2} \bigr]\varepsilon,
\end{eqnarray*}
with the second inequality following from Corollary~\ref{corollary 2.14.1}. The theorem is proved.
\end{pf}

%s7 #&#
%s7 ###
\section{Proof of Theorem \texorpdfstring{\protect\ref{theorem 5.19.1}}{2.1}}
\label{section 5.29.1}

According to Remark~\ref{remark 5.29.03}, in the proof of Theorem~\ref{theorem 5.19.1} we may assume that
$c^{\alpha\beta}(x)\geq\delta_{1}$.

First, we estimate the Lipschitz constant of $v$
on the boundary when $D\ne\bR^{d}$.

%
%le7.1 #&#
\begin{lemma}
\label{lemma 5.20.2}
Let $D$ be bounded and satisfy
the uniform exterior ball condition.
Let $x\in\bR^{d}$ and $y\notin D$.
Then there is a constant $N$ depending only on
$D$, $K_{0}$ and $\|g\|_{C^{2}(\bR^{d})}$, such that
\[
\bigl|v(x)-v(y)\bigr|\leq N|x-y|.
\]
\end{lemma}

\begin{pf}
If $x\notin D$, then
$|v(x)-v(y)|=|g(x)-g(y)|\leq N|x-y|$.
Therefore in the rest of the proof we assume that
$x\in D$. Then observe that by It\^o's formula we have
%
%
%e7.1 #&#
%e7.1 ###
\begin{equation}
\label{5.20.2} v(x)=g(x)+\infsup_{\bbeta\in\mathbb{B}\,  \alpha_{\cdot}\in\frA} E^{\alpha_{\cdot}\bbeta(\alpha_{\cdot})}_{x}
\int_{0}^{\tau} \bigl[L g(x_{t})
+f(x_{t}) \bigr]e^{-\phi_{t}} \,dt.
\end{equation}
It is well known that, in light
of the boundedness of $L^{\alpha\beta}g+f^{\alpha\beta}$
and $D$
and the uniform exterior ball condition,
the expectations in \eqref{5.20.2}
by magnitude are dominated by a constant
times $\dist(x,\partial D)\leq|x-y|$.
This proves the lemma since $v(y)=g(y)$ and $|g(x)-g(y)|
\leq N|x-y|$.
\end{pf}

\begin{pf*}{Proof of Theorem~\ref{theorem 5.19.1}}
In Section~\ref{section 5.20.1} take
\[
r\equiv1,\qquad p\equiv0,\qquad P\equiv I, \qquad\pi^{\alpha\beta}(x,y)=\mu \bigl[
\sigma^{\alpha\beta}(x) \bigr]^{*}(x-y),
\]
where the constant $\mu\geq1$ is chosen to be such
that \eqref{4.24.2}
with $\delta=1$ and
\[
(\sigma_{t},b_{t} ) (x,y)= (\sigma,b )^{\alpha_{t}\beta_{t}}(x )
\]
holds for all $\alpha_{\cdot}\in\frA$, $\beta_{\cdot}\in\frB$,
$x$ and $y$.
This is possible since $\sigma$ and $b$ are Lipschitz continuous,
and $a$ is uniformly nondegenerate. In
Section~\ref{section 5.20.1} we required $\pi^{\alpha\beta}(x,y)$
to be bounded and Lipschitz continuous
with respect to $x$. Since we will be
only concerned with its values for $|x-y|\leq1$, we can appropriately
modify the above
$\pi^{\alpha\beta}(x,y)$ for $|x-y|\geq1$
keeping the same notation.

Then for a unit $\xi\in\bR^{d}$,
$\varepsilon\geq0$, $\alpha_{\cdot}\in\frA$
and $\beta_{\cdot}\in\frB$
introduce $x^{\alpha_{\cdot}\beta_{\cdot}0}_{t}(
\varepsilon)$ as a unique solution of
\[
x_{t}= \varepsilon\xi+\int_{0}^{t}
\sigma^{\alpha_{s}\beta_{s}}(x_{s}) \,dw_{s} +\int
_{0}^{t} \bigl[ b^{\alpha_{s}\beta_{s}}(x_{s})-
\sigma^{\alpha_{s}\beta_{s}}(x_{s}) \pi^{\alpha_{s}\beta_{s}}(x_{s},y_{s})
\bigr] \,ds,
\]
where
\[
y_{s}=x^{\alpha_{\cdot}\beta_{\cdot}0}_{s}.
\]
Next introduce
\[
\phi_{t}^ {\alpha_{\cdot}\beta_{\cdot}0}(\varepsilon) =\int_{0}^{t}c^{\alpha_{s}\beta_{s}}
\bigl( x^{\alpha_{\cdot}\beta_{\cdot}0}_{s}( \varepsilon) \bigr) \,ds,
\]
and let $z_{t}^ {\alpha_{\cdot}\beta_{\cdot}0}(\varepsilon)$
be a unique solution of
\[
z_{t}=1+\int_{0}^{t}z_{s}
\bigl[\pi^{\alpha_{s},\beta_{s}} \bigl(x^{\alpha_{\cdot}\beta_{\cdot}0}_{s}(
\varepsilon),x^{\alpha_{\cdot}\beta_{\cdot}0}_{s}(0) \bigr) \bigr]^{*}
\,dw_{s}.
\]

Keeping in mind that $\mu$ is already fixed, set $\delta_{1}:=
\varepsilon_{1}=1$, take $\lambda$ from Lemma~\ref{lemma 4.7.3},
fix $\varepsilon\in(0,1]$
and introduce
\begin{eqnarray*}
\tau^{\alpha_{\cdot}\beta_{\cdot}0}_{\varepsilon} &=&\inf \bigl\{t\geq0\dvtx x^{\alpha_{\cdot}\beta_{\cdot}0}_{t}(
\varepsilon) \notin D \bigr\},
\\
\gamma_{\varepsilon}^{\alpha_{\cdot}\beta_{\cdot}0} &=&\inf \bigl\{t\geq0\dvtx
z^{\alpha_{\cdot}\beta_{\cdot}0}_{t}( \varepsilon)e^{\phi_{t}^{\alpha_{\cdot}\beta_{\cdot}0}
(0)-\phi_{t}^{\alpha_{\cdot}\beta_{\cdot}0}
(\varepsilon)}\notin(1/2,2) \bigr
\},
\\
\kappa_{\varepsilon}^{\alpha_{\cdot}\beta_{\cdot}0} &=&\inf \bigl\{t\geq0\dvtx
\bigl|x^{\alpha_{\cdot}\beta_{\cdot}0}_{t}( \varepsilon)-
x^{\alpha_{\cdot}\beta_{\cdot}0}_{t}(
0)\bigr|\geq\lambda \bigr\},
\\
\gamma^{\alpha_{\cdot}\beta_{\cdot}0} &=&
\tau^{\alpha_{\cdot}\beta_{\cdot}0}_{\varepsilon} \wedge
\tau^{\alpha_{\cdot}\beta_{\cdot}0}_{0} \wedge\kappa_{\varepsilon}^{\alpha_{\cdot}\beta_{\cdot}0}
\wedge\gamma_{\varepsilon}^{\alpha_{\cdot}\beta_{\cdot}0}.
\end{eqnarray*}

By Theorem~\ref{theorem 5.1.1},
%
%
%e7.2 #&#
%e7.2 ###
\begin{eqnarray}
\label{5.19.4}&& v(\varepsilon\xi) =\infsup_{\bbeta\in\mathbb{B}\,  \alpha_{\cdot}\in\frA} E_{0}^{\alpha_{\cdot}\bbeta(\alpha_{\cdot})}
\biggl[ \int_{0}^{\gamma}z_{t}(\varepsilon)
f \bigl(x_{t}(\varepsilon) \bigr) e^{- \phi_{t}(\varepsilon) } \,dt
\nonumber
\\[-8pt]
\\[-8pt]
\nonumber
&&\hspace*{104pt}\qquad{}+
z_{\gamma}(\varepsilon)v \bigl(x_{\gamma}(\varepsilon) \bigr)
e^{- \phi_{\gamma}(\varepsilon) } \biggr].
\end{eqnarray}

Next we fix $\alpha_{\cdot}\in\frA$ and $\beta_{\cdot}\in\frB$,
and
in Section~\ref{section 5.20.2} use the functions
\[
(\sigma_{t},b_{t},c_{t},f_{t})
(x,y)= (\sigma,b,c,f)^{\alpha_{t}\beta_{t}}(x).
\]
Observe that in the expectation
\[
E_{0}^{\alpha_{\cdot}\beta_{\cdot}} \biggl[ \int_{0}^{\gamma}z_{t}(
\varepsilon) f \bigl(x_{t}(\varepsilon) \bigr) e^{- \phi_{t}(\varepsilon) } \,dt +
z_{\gamma}(\varepsilon)v \bigl(x_{\gamma}(\varepsilon) \bigr)
e^{- \phi_{\gamma}(\varepsilon) } \biggr],
\]
one can replace $x^{\alpha_{\cdot}\beta_{\cdot}0}_{s}(
\varepsilon)$ with $x^{\varepsilon}_{t}$ since both satisfy
the same equation on $[0,\gamma^{\alpha_{\cdot}\beta_{\cdot}0}]$,
and by Theorem~\ref{theorem 4.6.1}
we get that
%
%e7.3 ###
\begin{eqnarray}\label{5.20.4}
&&\biggl| E_{0}^{\alpha_{\cdot}\beta_{\cdot}} \biggl[ \int_{0}^{\gamma}z_{t}(
\varepsilon) f \bigl(x_{t}(\varepsilon) \bigr) e^{- \phi_{t}(\varepsilon) } \,dt +
z_{\gamma}(\varepsilon)v \bigl(x_{\gamma}(\varepsilon) \bigr)
e^{- \phi_{\gamma}(\varepsilon) } \biggr]
\nonumber\\
&&\hspace*{61pt}\qquad{}-E_{0}^{\alpha_{\cdot}\beta_{\cdot}} \biggl[ \int_{0}^{\gamma}
f(x_{t} ) e^{- \phi_{t} } \,dt + v (x_{\gamma} )
e^{- \phi_{\gamma} } \biggr] \biggr|
\\
 &&\qquad\leq N\varepsilon +E_{0}^{\alpha_{\cdot}\beta_{\cdot}} \bigl|v \bigl(x
_{\gamma}(\varepsilon) \bigr)- v \bigl(x _{\gamma}(0)
\bigr)\bigr|e^{-\phi
_{\gamma}} I_{\gamma<\gamma_{\varepsilon}
\wedge\kappa_{\varepsilon} }.\nonumber
\end{eqnarray}

If $t=\gamma^{\alpha_{\cdot}\beta_{\cdot}0}
<\gamma^{\alpha_{\cdot}\beta_{\cdot}0}_{\varepsilon}\wedge
\kappa^{\alpha_{\cdot}\beta_{\cdot}0}_{\varepsilon}$,
then ($D\ne\bR^{d}$ and) at least one of $x^{\alpha_{\cdot}\beta
_{\cdot}0}_{t}(
\varepsilon)$ and $x^{\alpha_{\cdot}\beta_{\cdot}0}_{t}(0)$
is outside $D$, and by Lemma~\ref{lemma 5.20.2} we obtain
\begin{eqnarray*}
&&E_{0}^{\alpha_{\cdot}\beta_{\cdot}} \bigl|v \bigl(x _{\gamma}(\varepsilon) \bigr)-
v \bigl(x _{\gamma}(0) \bigr)\bigr|e^{-\phi
_{\gamma}} I_{\gamma<\gamma_{\varepsilon}
\wedge\kappa_{\varepsilon} } \\
&&\qquad\leq N
E_{0}^{\alpha_{\cdot}\beta_{\cdot}} \bigl|x _{\gamma}(\varepsilon) - x
_{\gamma}(0) \bigr|e^{-\phi
_{\gamma}} I_{\gamma<\gamma_{\varepsilon}
\wedge\kappa_{\varepsilon} }
\\
&&\qquad=N\varepsilon E_{0}^{\alpha_{\cdot}\beta_{\cdot}} \bigl|\xi_{\gamma}(
\varepsilon)\bigr|e^{-\phi
_{\gamma}} I_{\gamma<\gamma_{\varepsilon}
\wedge\kappa_{\varepsilon} }\leq N\varepsilon
E_{0}^{\alpha_{\cdot}\beta_{\cdot}}\sup_{t<
\kappa_{\varepsilon} }\bigl| \xi_{t}(
\varepsilon)\bigr|e^{-\phi
_{t}},
\end{eqnarray*}
where $\varepsilon\xi^{\alpha_{\cdot}\beta_{\cdot}0}_{t}(
\varepsilon)=x^{\alpha_{\cdot}\beta_{\cdot}0}_{t}(
\varepsilon)-x^{\alpha_{\cdot}\beta_{\cdot}0}_{t}(
0)$. By using Lemma~\ref{lemma 2.3.1}, equation~\eqref{5.20.4},
and the fact that $\alpha_{\cdot}$ and $\beta_{\cdot}$
in the above argument are arbitrary, we see that
$|v(\varepsilon\xi)-v(0)|\leq N\varepsilon$.
Similarly one proves that
$|v(x+\varepsilon\xi)-v(x)|\leq N\varepsilon$ for
any $x$, which is what
we need. The theorem is proved.
\end{pf*}

%s8 #&#
%s8 ###
\section{Proof of Theorem \texorpdfstring{\protect\ref{theorem 6.4.1}}{2.2}}
\label{section 6.4.1}

In contrast with Section~\ref{section 5.29.1}, where we used
$\delta=1$, here $\delta$ will be chosen large.
We begin with the following.
%
%
%le8.1 #&#
\begin{lemma}
\label{lemma 6.4.1}
Let $D$ be a bounded domain satisfying the uniform
exterior ball condition, and let $\|g\|_{C^{2}(\bR^{d})}
<\infty$. For $R\in(0,1]$
let $B_{R}=\{x\dvtx |x|\leq R\}$.
Assume that for an $R$ we have
$B_{R}\subset D$ and denote by $L_{R}$
the Lipschitz constant of $v$ in $B_{R}$
(finite by Theorem~\ref{theorem 5.19.1}).
Finally assume that $|v|\leq K_{0}$ in $B_{R}$.

Then for any $\delta\geq K_{1}^{2}+4K^{2}_{0}+2$ we have
%
%
%e8.1 #&#
%e8.1 ###
\begin{equation}
\label{6.4.3} \nlimsup_{x\to0}\frac{|v(x)-v(0)|}{|x|} \leq N\delta
R^{-1}+Ne^{-\nu\sqrt{\delta}}L_{R},
\end{equation}
where $N$ and $\nu>0$ depend only on $d$,
$K_{0}$, $K_{1}$ and $\delta_{0}$.
\end{lemma}
\begin{pf} First suppose that $R=1$.
Observe that by the dynamic programming principle
%
%
%e8.2 #&#
%e8.2 ###
\begin{equation}
\label{6.5.1} v(x) =\infsup_{\bbeta\in\mathbb{B}\,  \alpha_{\cdot}\in\frA} E_{x}^{\alpha_{\cdot}\bbeta(\alpha_{\cdot})}
\biggl[ \int_{0}^{\tau_{1}} f(x_{t} )
e^{- \phi_{t} } \,dt + v \bigl(x_{\tau_{1}}(\varepsilon)
\bigr)e^{- \phi_{\tau_{1}} } \biggr],
\end{equation}
where $\tau_{1}^{\alpha_{\cdot}\beta_{\cdot}x}$
is the first exit time of
$x^{\alpha_{\cdot}\beta_{\cdot}x}_{t}$ from $B_{1}$.

Remark~\ref{remark 5.29.03} allows us to rewrite
\eqref{6.5.1} by using a global barrier for $B_{1}$
for a slightly modified $v$. Obviously, if we can prove
\eqref{6.4.3} with $R=1$ for such modification,
then we will have it also for the original function.
Hence, concentrating on~\eqref{6.5.1}
\emph{and the case $R=1$}, without losing
generality we may assume that $c^{\alpha\beta}\geq1$.

Set $\mu= \delta_{0}^{-1}
\delta+N_{0}$, where $N_{0}$ depending only on $K_{1}$,
$\delta_{0}$, and $d$ is chosen in such a way
that \eqref{4.24.2} is
satisfied with
\[
(\sigma_{t},b_{t} ) (x,y)= (\sigma,b )^{\alpha_{t}\beta_{t}}(x )
\]
for all $\alpha_{\cdot}\in\frA$, $\beta_{\cdot}\in\frB$,
$x$, $y$ and $\delta>0$.

We use the notation from
the proof of Theorem~\ref{theorem 5.19.1} in
Section~\ref{section 5.29.1} and write \eqref{5.19.4}
with
\[
\gamma^{\alpha_{\cdot}\beta_{\cdot}0} =\tau^{\alpha_{\cdot}\beta_{\cdot}0}_{1}(\varepsilon) \wedge
\tau^{\alpha_{\cdot}\beta_{\cdot}0}_{1}(0) \wedge\kappa_{\varepsilon}^{\alpha_{\cdot}\beta_{\cdot}0}
\wedge\gamma_{\varepsilon}^{\alpha_{\cdot}\beta_{\cdot}0},
\]
where $\tau^{\alpha_{\cdot}\beta_{\cdot}0}_{1}(\varepsilon)$
is the first exit time of
$x^{\alpha_{\cdot}\beta_{\cdot}0}_{t}(
\varepsilon)$ from $B_{1}$.

As in the proof of Theorem~\ref{theorem 5.19.1},
by Theorem~\ref{theorem 4.6.1} (with $\tau=
\gamma^{\alpha_{\cdot}\beta_{\cdot}0}$ there),
we get that
(recall that $M=2$ and $\mu$ is of order
$\delta$ if $\delta\geq1$)
%
%
%e8.3 #&#
%e8.3 ###
\begin{equation}
\label{6.4.5} \bigl|v(\varepsilon\xi)-v(0)\bigr|\leq N\delta\varepsilon+S_{\varepsilon},
\end{equation}
where $N$ depends only on
$K_{0}$, $K_{1}$ and $\delta_{0}$ (recall that
$\delta_{1}=1$) and
\begin{eqnarray*}
S_{\varepsilon}&:=& \sup_{\alpha_{\cdot},\beta_{\cdot}} E^{\alpha_{\cdot}\beta_{\cdot}}_{0}
\bigl|v \bigl(x_{\gamma}(\varepsilon) \bigr) -v \bigl(x_{\gamma}(0)
\bigr)\bigr|e^{-\phi_{\gamma}}I_{
\tau_{1}(\varepsilon)\wedge\tau_{1}(0)<
\gamma_{\varepsilon}\wedge\kappa_{\varepsilon}}
\\
&\leq&\varepsilon L_{1} \sup_{\alpha_{\cdot},\beta_{\cdot}}
E^{\alpha_{\cdot}\beta_{\cdot}}_{0} \bigl|\xi_{
 \tau_{1}(\varepsilon)\wedge\tau_{1}(0)}(\varepsilon)\bigr|
e^{-\phi_{\tau_{1}(\varepsilon)\wedge\tau_{1}(0)}}I_{
\tau_{1}(\varepsilon)\wedge\tau_{1}(0)<
\kappa_{\varepsilon}}.
\end{eqnarray*}

Observe that for any $T>0$ by Lemma~\ref{lemma 2.3.1}
and Remark~\ref{remark 6.4.1} ($\delta\geq K_{1}^{2}$),
\begin{eqnarray*}
&&E^{\alpha_{\cdot}\beta_{\cdot}}_{0} \bigl|\xi_{ \tau_{1}(\varepsilon)}
(\varepsilon)\bigr|
e^{-\phi_{ \tau_{1}(\varepsilon)}}I_{
\tau_{1}(\varepsilon)<
\kappa_{\varepsilon}}\\
&&\qquad= E^{\alpha_{\cdot}\beta_{\cdot}}_{\varepsilon}\bigl |
\xi_{ \tau_{1}(\varepsilon)}(\varepsilon)\bigr| e^{-\phi_{ \tau_{1}(\varepsilon)}}I_{
\tau_{1}(\varepsilon)<
\kappa_{\varepsilon}\wedge T}
\\
&&\qquad\quad{}+E^{\alpha_{\cdot}\beta_{\cdot}}_{0} \bigl|\xi_{ \tau_{1}(\varepsilon)}
(\varepsilon)\bigr|
e^{-\phi_{ \tau_{1}(\varepsilon)}}I_{
\tau_{1}(\varepsilon)<
\kappa_{\varepsilon}}I_{
\tau_{1}(\varepsilon)\geq T}
\\
&&\qquad\leq \bigl(E^{\alpha_{\cdot}\beta_{\cdot}}_{0}I_{\tau
_{1}(\varepsilon)
<T}
\bigr)^{1/2} +e^{-\delta T/2}E^{\alpha_{\cdot}\beta_{\cdot}}_{0} \sup
_{t<\kappa_{\varepsilon}}\bigl|\xi_{t}(\varepsilon)\bigr| e^{-\phi_{t}+\delta t/2}
\\
&&\qquad\leq Ne^{-\delta T/2}+ \bigl(E^{\alpha_{\cdot}\beta_{\cdot}}_{0} I_{\tau_{1}(\varepsilon)
<T}
\bigr)^{1/2}.
\end{eqnarray*}
Similarly,
\[
E^{\alpha_{\cdot}\beta_{\cdot}}_{0} \bigl|\xi_{ \tau_{1}(0)}(\varepsilon)\bigr|
e^{-\phi_{ \tau_{1}(0)}}I_{
\tau_{1}(0)<
\kappa_{\varepsilon}} \leq Ne^{-\delta T/2}+ \bigl(E^{\alpha_{\cdot}\beta_{\cdot}}_{0}
I_{\tau_{1}(0)
<T} \bigr)^{1/2}.
\]

One knows that if the starting point of a diffusion process
with coefficients bounded by $K_{0}$ is in the ball
of radius $\varepsilon<1/2$, then the probability
that the process will exit from $B_{1}$ before time $T$
is less than $N\exp(-\nu/T)$ if $K_{0}T\leq1/2$, where $N$ and
$\nu$ depend only
on $K_{0}$ and $d$. This result is easily obtained by using
the McKeen estimate (see, e.g., Corollary
IV.2.9 of \cite{Intro}) for each coordinate
of the process from which one subtracts the drift term.
Hence (with another $\nu$)
\[
S_{\varepsilon}\leq\varepsilon L_{1} \bigl( Ne^{-\delta T/2}+Ne^{-\nu/T}
\bigr).
\]
For $T= \delta^{-1/2}$ (so that
$K_{0}T\leq1/2$ since $\delta\geq4K^{2}_{0}$)
we get that
(yet with another $\nu$)
\[
S_{\varepsilon}\leq\varepsilon L_{1} Ne^{-\nu\sqrt{\delta}},
\]
and the result follows in case $R=1$.

Once \eqref{6.4.3} is proved for $R=1$,
for $R\in(0,1)$ it
follows by using dilations (see Remark~2.5 of \cite{Kr_2}),
which allow us to keep the
constants $\delta_{0},K_{0}$ and $K_{1}$ (actually,
after dilations the constant $K_{1}$
can be taken even smaller then the original one).
The lemma is proved.
\end{pf}

\begin{pf*}{Proof of Theorem~\ref{theorem 6.4.1}}
First suppose that $\|g\|_{C^{2}(\bR^{d})}
<\infty$ and
that for an $R_{0}>0$ we have $B_{2R_{0}}
\subset D$. Estimate \eqref{6.4.3} can be applied to any point
rather than only $0$, and it shows that
for any $R'<R''\leq2R_{0}$ and $\delta\geq
K_{1}^{2}+4K^{2}_{0}+2$
we have
\[
L_{R'}\leq N\delta/ \bigl(R''-R'
\bigr)+N_{1}e^{-\nu\sqrt{\delta}} L_{R''}.
\]
We apply this inequality to $R'=R_{n}$ and $R''=R_{n+1}$,
where $R_{n}$, $n\geq1$, are defined by
\[
R_{n}=R_{0}+R_{0}\sum
_{i=1}^{n}\frac{\chi}{i^{2}},
\]
and $\chi$ is such that $R_{n}\to2R_{0}$ as $n\to\infty$.
We also take and fix $\delta\geq K_{1}^{2}+4K^{2}_{0}+2$
so large that $N_{1}e^{-\nu\sqrt{\delta}}\leq1/2$. Then for
a constant $N_{0}$ depending only on $\delta_{0},K_{0},K_{1}$
and $d$ and all $n\geq0$, we get that
\begin{eqnarray*}
L_{R_{n}}&\leq& N_{0}R_{0}^{-1}(n+1)^{2}+2^{-1}L_{R_{n+1}},
\\
2^{-n}L_{R_{n}}&\leq&2^{-n}N_{0}R_{0}^{-1}
(n+1)^{2}+2^{-(n+1)}L_{R_{n+1}},
\\
\sum_{n=0}^{\infty} 2^{-n}L_{R_{n}}
&\leq& N_{0}R_{0}^{-1}\sum
_{n=0}^{\infty} 2^{-n}(n+1)^{2}+
\sum_{n=0}^{\infty}2^{-(n+1)}L_{R_{n+1}}
\end{eqnarray*}
and $L_{R_{0}} \leq N_{0}IR_{0}^{-1}$, where
\[
I=2\sum_{n=1}^{\infty}2^{-n}n^{2}.
\]
One can do the same estimate for any ball
inside $D$ not necessarily centered at the origin,
and this yields the desired result in case
$\|g\|_{C^{2}(\bR^{d})}
<\infty$. In the general case where $g$ is only continuous,
it suffices to use appropriate approximations of
it by smooth functions. The theorem is proved.
\end{pf*}

%s9 #&#
%s9 ###
\section{Proof of Theorem \texorpdfstring{\protect\ref{theorem 5.22.1}}{2.3}}
\label{section 6.4.2}

First of all we point out that the assertion of Lemma~\ref{lemma 5.20.2} continues to hold true
with only one difference that $N$ depends only on
$K_{0}$, $G$, $d$ and $\|g\|_{C^{2}(\bR^{d})}$. The proof
remains the same with It\^o's formula showing
that the expectations in \eqref{5.20.2} are bounded
by $NG(x)$. The remaining arguments follow the ones
from Section~\ref{section 5.29.1}
almost word for word.

In Section~\ref{section 5.20.1} for $|x-y|\leq1$ take
\[
\pi^{\alpha\beta}(x,y)= \mu \bigl[\sigma^{\alpha\beta}(y) \bigr]^{*}(x-y)
\]
and extend it appropriately for $|x-y|>1$.

Then for a unit $\xi\in\bR^{d}$,
$\varepsilon\geq0$, $\alpha_{\cdot}\in\frA$,
and $\beta_{\cdot}\in\frB$
introduce $x^{\alpha_{\cdot}\beta_{\cdot}0}_{t}(\varepsilon)$
as a unique solution of equation \eqref{5.7.2}
with initial condition $\varepsilon\xi$ and
\[
y_{s}=x^{\alpha_{\cdot}\beta_{\cdot}0}_{s}.
\]
Observe that $x^{\alpha_{\cdot}\beta_{\cdot}0}_{t}(0)
=x^{\alpha_{\cdot}\beta_{\cdot}0}_{t}$. Then
define $z_{t}^ {\alpha_{\cdot}\beta_{\cdot}0}(\varepsilon)$,
$\tau^{\alpha_{\cdot}\beta_{\cdot}0}_{\varepsilon}$,
$\gamma_{\varepsilon}^{\alpha_{\cdot}\beta_{\cdot}0} $,
$\kappa_{\varepsilon}^{\alpha_{\cdot}\beta_{\cdot}0} $
and $\gamma^{\alpha_{\cdot}
\beta_{\cdot}0}$ in the same way as in Section~\ref{section 5.29.1},
and use
Theorem~\ref{theorem 5.1.1} to get that
\begin{eqnarray*}
&&v(\varepsilon\xi) =\infsup_{\bbeta\in\mathbb{B}\,  \alpha_{\cdot}\in\frA} E_{0}^{\alpha_{\cdot}\bbeta(\alpha_{\cdot})}
\biggl[
z_{\gamma}(\varepsilon)v \bigl(x_{\gamma}(\varepsilon) \bigr)
e^{- \phi_{\gamma}(\varepsilon) }
\\
&&\hspace*{90pt}\qquad{}+\int_{0}^{\gamma}z_{t}(\varepsilon) \hat
f \bigl(x_{t}(\varepsilon),x_{t}(0) \bigr) e^{- \phi_{t}(\varepsilon) }
\,dt \biggr],
\end{eqnarray*}
where
\[
\phi_{t}^ {\alpha_{\cdot}\beta_{\cdot}0}(\varepsilon) =\int_{0}^{t}
\hat c^{\alpha_{s}\beta_{s}} \bigl( x^{\alpha_{\cdot}\beta_{\cdot}0}_{s}(
\varepsilon),x^{\alpha_{\cdot}\beta_{\cdot}0}_{s}(0) \bigr) \,ds.
\]

Fix $\alpha_{\cdot}\in\frA$ and $\beta_{\cdot}\in\frB$,
and
in Section~\ref{section 5.20.2} use the functions
\[
(\sigma_{t},b_{t},c_{t},f_{t})
(x,y)= (\hat\sigma,\hat b,\hat c,\hat f)^{\alpha_{t}\beta_{t}}(x,y).
\]
Observe that Assumption~\ref{assumption 4.24.1}
is satisfied owing to Assumption~\ref{assumption 5.29.2}.

Furthermore,
for $t\leq\gamma^{\alpha_{\cdot}\beta_{\cdot}}$
the processes $x^{\varepsilon}_{t}$ and $y_{t}$
coincide with $x^{\alpha_{\cdot}\beta_{\cdot}0}_{t}
(\varepsilon)$ and $x^{\alpha_{\cdot}\beta_{\cdot}0}_{t}
(0)$, respectively, since they satisfy the same equations,
respectively. It follows that in the expectation
\[
E_{0}^{\alpha_{\cdot}\beta_{\cdot}} \biggl[ \int_{0}^{\gamma}z_{t}(
\varepsilon) f \bigl(x_{t}(\varepsilon),x_{t}(0) \bigr)
e^{- \phi_{t}(\varepsilon) } \,dt + z_{\gamma}(\varepsilon)v \bigl(x_{\gamma}(
\varepsilon) \bigr) e^{- \phi_{\gamma}(\varepsilon) } \biggr],
\]
one can replace $x^{\alpha_{\cdot}\beta_{\cdot}0}_{s}(
\varepsilon)$ with $x^{\varepsilon}_{t}$,
and by Theorem~\ref{theorem 4.6.1}
we get that
%
%e9.1 ###
\begin{eqnarray}\label{5.20.04}
&&\biggl| E_{0}^{\alpha_{\cdot}\beta_{\cdot}} \biggl[ \int_{0}^{\gamma}z_{t}(
\varepsilon) f \bigl(x_{t}(\varepsilon),x_{t}(0) \bigr)
e^{- \phi_{t}(\varepsilon) } \,dt + z_{\gamma}(\varepsilon)v \bigl(x_{\gamma}(
\varepsilon) \bigr) e^{- \phi_{\gamma}(\varepsilon) } \biggr]
\nonumber\\
&&\hspace*{109pt}{}-E_{0}^{\alpha_{\cdot}\beta_{\cdot}} \biggl[ \int_{0}^{\gamma}
f(x_{t} ) e^{- \phi_{t} } \,dt + v (x_{\gamma} )
e^{- \phi_{\gamma} } \biggr]\biggr|
\\
&&\qquad\leq N\varepsilon +E_{0}^{\alpha_{\cdot}\beta_{\cdot}}\bigl |v \bigl(x
_{\gamma}(\varepsilon) \bigr)- v \bigl(x _{\gamma}(0)
\bigr)\bigr|e^{-\phi
_{\gamma}} I_{\gamma<\gamma_{\varepsilon}
\wedge\kappa_{\varepsilon} }.\nonumber
\end{eqnarray}

If $t=\gamma^{\alpha_{\cdot}\beta_{\cdot}0}
<\gamma^{\alpha_{\cdot}\beta_{\cdot}0}_{\varepsilon}\wedge
\kappa^{\alpha_{\cdot}\beta_{\cdot}0}_{\varepsilon}$,
then at least one of $x^{\alpha_{\cdot}\beta_{\cdot}0}_{t}(
\varepsilon)$ and $x^{\alpha_{\cdot}\beta_{\cdot}0}_{t}(0)$
is outside $D$, and by Lemma~\ref{lemma 5.20.2} we obtain
\begin{eqnarray*}
&&E_{0}^{\alpha_{\cdot}\beta_{\cdot}}\bigl |v \bigl(x _{\gamma}(\varepsilon) \bigr)-
v \bigl(x _{\gamma}(0) \bigr)\bigr|e^{-\phi
_{\gamma}} I_{\gamma<\gamma_{\varepsilon}
\wedge\kappa_{\varepsilon} } \\
&&\qquad\leq N
E_{0}^{\alpha_{\cdot}\beta_{\cdot}} \bigl|x _{\gamma}(\varepsilon) - x
_{\gamma}(0) \bigr|e^{-\phi
_{\gamma}} I_{\gamma<\gamma_{\varepsilon}
\wedge\kappa_{\varepsilon} }
\\
&&\qquad=\varepsilon E_{0}^{\alpha_{\cdot}\beta_{\cdot}} \bigl|\xi_{\gamma}(
\varepsilon)\bigr|e^{-\phi
_{\gamma}} I_{\gamma<\gamma_{\varepsilon}
\wedge\kappa_{\varepsilon} }\leq \varepsilon
E_{0}^{\alpha_{\cdot}\beta_{\cdot}}\sup_{t<
\kappa_{\varepsilon} }\bigl| \xi_{t}(
\varepsilon)\bigr|e^{-\phi
_{t}},
\end{eqnarray*}
where $\varepsilon\xi^{\alpha_{\cdot}\beta_{\cdot}0}_{t}(
\varepsilon)=x^{\alpha_{\cdot}\beta_{\cdot}0}_{t}(
\varepsilon)-x^{\alpha_{\cdot}\beta_{\cdot}0}_{t}(
0)$. By using Lemma~\ref{lemma 2.3.1}, \eqref{5.20.04}
and the fact that $\alpha_{\cdot}$ and $\beta_{\cdot}$
in the above argument are arbitrary, we see that
$|v(\varepsilon\xi)-v(0)|\leq N\varepsilon$.
Similarly one proves that
$|v(x+\varepsilon\xi)-v(x)|\leq N\varepsilon$ for
any $x$, which is what
we need. The theorem is proved.

%s10 #&#
%s10 ###
\section{Proof of Theorem \texorpdfstring{\protect\ref{theorem 6.3.1}}{2.4}}
\label{section 6.8.1}

Obviously $v\leq v_{K}$. To estimate $v_{K}-v$ from above,
define
\[
d_{K}=\sup_{\bR^{d}}(v_{K}-v),\qquad \lambda=
\sup_{\alpha,\beta,x} c^{\alpha\beta}(x).
\]
By the dynamic programming principle (see Theorem~3.1 in
\cite{Kr_2}),
\[
v_{K}(x)=\infsup_{\bbeta\in\hat{B}\,  \alpha_{\cdot}\in\hat
{\frA}} E_{x}^{\alpha_{\cdot}\bbeta(\alpha_{\cdot})}
\biggl[v_{K}(x_{1})e^{ -\lambda} +\int
_{0}^{1} \bigl\{f_{K}+(\lambda-c )
v_{K} \bigr\}(x_{t}) e^{ -\lambda t} \,dt \biggr].
\]
Observe that
\[
e^{-\lambda}+\int_{0}^{1} \bigl[
\lambda-c^{\alpha_{t}\beta_{t}} \bigl(x^{\alpha_{\cdot}\beta_{\cdot}x}_{t} \bigr)
\bigr]e^{-\lambda t} \,dt \leq e^{-\lambda}+\int_{0}^{1}(
\lambda-\delta_{1}) e^{-\lambda t} \,dt =:\kappa<1.
\]
Hence,
\[
v_{K}(x)\leq\infsup_{\bbeta\in\hat{B}\,  \alpha_{\cdot}\in\hat
{\frA}} E_{x}^{\alpha_{\cdot}\bbeta(\alpha_{\cdot})}
\biggl[v (x_{1})e^{ -\lambda} +\int_{0}^{1}
\bigl\{f_{K}+(\lambda-c ) v \bigr\}(x_{t})
e^{ -\lambda t} \,dt \biggr]+\kappa d_{K}.
\]

Now take a sequence $x^{n}$ maximizing $v_{K}-v$, and
take $\bbeta^{n}\in\mathbb{B}$ such that
%
%
%e10.1 #&#
%e10.1 ###
\begin{eqnarray}
\label{6.3.3} v \bigl(x^{n} \bigr)&\geq&\sup_{\alpha_{\cdot}\in\frA}
E^{\alpha_{\cdot}\bbeta^{n}(\alpha_{\cdot})}_{x^{n}} \biggl[\int_{0}^{1}
\bigl(f+(\lambda-c)v \bigr) (x_{t})e^{-\lambda t} \,dt
+e^{-\lambda}v(x_{1}) \biggr]
\nonumber
\\[-8pt]
\\[-8pt]
\nonumber
&&{}-1/n.
\end{eqnarray}
Also define $\pi\alpha=\alpha$
if $\alpha\in A_{1}$ and $\pi\alpha=\alpha^{*}$
if $\alpha\in A_{1}$, where $\alpha^{*}$ is a fixed
element of $A_{1}$, and find $\alpha_{\cdot}^{n}\in\hat{\frA}$
such that
%
%e10.2 ###
\begin{eqnarray}\label{6.3.4}
v_{K} \bigl(x^{n} \bigr)&\leq& E_{x^{n}}^{\alpha^{n}_{\cdot}\bbeta^{n}(\pi\alpha^{n}_{\cdot})}
\biggl[v (x_{1})e^{ -\lambda} +\int_{0}^{1}
\bigl\{f_{K}+(\lambda-c ) v \bigr\}(x_{t})
e^{ -\lambda t} \,dt \biggr]\nonumber\\
&&{}+\kappa d_{K}+1/n
\\
\nonumber& =&E_{x^{n}}^{\alpha^{n}_{\cdot}\bbeta^{n}(\pi\alpha^{n}_{\cdot})} \biggl[v (x_{1})e^{ -\lambda}
+\int_{0}^{1} \bigl\{f +(\lambda-c ) v \bigr
\}(x_{t}) e^{ -\lambda t} \,dt \biggr]
\\
&&{}-KR_{n} +\kappa d_{K}+1/n,\nonumber
\end{eqnarray}
where
\[
R_{n}=E\int_{0}^{1}e^{-\lambda t}I_{\alpha^{n}_{t}\in A_{2}}
\,dt.
\]

By Lemma~5.3 of \cite{Kr_2} for any $\alpha_{\cdot}
\in\hat\frA$, $\beta_{\cdot}\in\frB$ and $x\in\bR^{d}$, we have
\[
E\sup_{t\leq1}\bigl|x_{t}^{\pi\alpha_{\cdot}\beta_{\cdot}x}
-x_{t}^{ \alpha_{\cdot}\beta_{\cdot}x}\bigr|\leq N \biggl(E^{\alpha_{\cdot}\beta_{\cdot}}_{x}
\int_{0}^{1} e^{-t}I_{\alpha^{n}_{t}\in A_{2}}
\,dt \biggr)^{1/2},
\]
where the constant $N$ depends only on $K_{0}$,
$K_{1}$ and $d$.
We use this, and
since $c,f,v$ are Lipschitz continuous, we get from
\eqref{6.3.4} and \eqref{6.3.3},
\begin{eqnarray*}
&&v_{K} \bigl(x^{n} \bigr)+(K-N_{0})R_{n}
\\
&&\qquad\leq E_{x^{n}}^{\pi\alpha^{n}_{\cdot}\bbeta^{n}(\pi\alpha^{n}_{\cdot})} \biggl[v (x_{1})e^{ -\lambda}
+\int_{0}^{1} \bigl\{f +(\lambda-c ) v \bigr
\}(x_{t}) e^{ -\lambda t} \,dt \biggr]\\
&&\qquad\quad{}+\kappa d_{K}+1/n+NR^{1/2}_{n}\\
&&\qquad\leq v
\bigl(x^{n} \bigr)+\kappa d_{K}+2/n+NR^{1/2}_{n},
\end{eqnarray*}
where the constant $N_{0}$ depends only on the supremums of $c$,
$v$ and $f$.
Hence
%
%
%e10.3 #&#
%e10.3 ###
\begin{equation}
\label{6.3.6} v_{K} \bigl(x^{n} \bigr)-v
\bigl(x^{n} \bigr)-\kappa d_{K}+(K-N_{0})R_{n}
\leq2/n+NR^{1/2}_{n}.
\end{equation}
When $n$ is large enough, $v_{K}(x^{n})-v(x^{n})-\kappa d_{K}\geq0$
because of the way we chose $x^{n}$ and the fact that $\kappa<1$.
It follows that for $n$ large enough,
\[
(K-N_{0})R_{n}\le2/n+NR^{1/2}_{n},
\]
which for $K\geq2N_{0}+1$ implies that $K
R_{n}\le4/n+NR^{1/2}_{n}$, so that, if $KR_{n}\geq
8/n$, then $KR_{n}\leq NR_{n}^{1/2}$ and
$R_{n}\leq N/K^{2}$. Thus
\[
R_{n}\leq8/(nK)+N/K^{2},
\]
which after coming back to \eqref{6.3.6} finally yields
\begin{eqnarray*}
v_{K} \bigl(x^{n} \bigr)-v \bigl(x^{n}
\bigr)- \kappa d_{K} &\leq&2/n+N/\sqrt{n}+N/K,
\\
(1-\kappa)d_{K}&=&\lim_{n\to\infty} \bigl[v_{K}
\bigl(x^{n} \bigr)-v \bigl(x^{n} \bigr) \bigr]-\kappa
d_{K} \leq N/K,
\end{eqnarray*}
and the theorem is proved.
% imsref loaded by akundreckaite, 2014-01-16 08:28:22
% imsref loaded by akundreckaite, 2014-01-16 08:59:52
% imsref loaded by akundreckaite, 2014-01-20 16:23:05

% zodis "Acknowledgments" paliekamas pagal autoriu

%suskaldyti doi

\printaddresses

\end{document}